\documentclass{amsart}

\newif\ifboxes
\boxestrue 

\usepackage{amsaddr}
\usepackage{wrapfig}
\usepackage{newpxtext,newpxmath}

\usepackage{tikz-cd}

\usepackage[dvipsnames,table]{xcolor}

\definecolor{burntsienna}{rgb}{0.91, 0.45, 0.32}

\usepackage[backend=biber, 
style=alphabetic, 
maxcitenames=5, 
maxbibnames=9,
sorting=nyt]{biblatex}
 \theoremstyle{definition}
\newtheorem*{remarks*}{Remarks}

\usepackage[margin=3.2cm]{geometry}

\renewcommand{\S}{\mathcal{S}}

\newcommand{\ubull}[1]{{(#1)_<}}
\newcommand{\obull}[1]{(#1)_>}

\newcommand{\un}{n_<}

\newcommand{\uk}{k_<}
\newcommand{\on}{n_>}
\newcommand{\ok}{k_>}

\newcommand{\wtm}{\widetilde{m}}

\newcommand{\ff}[2]{\left\lfloor\frac{#1}{#2}\right\rfloor}

\renewcommand{\mod}{\ \mathrm{mod}\,}

\usepackage{keytheorems}
\usepackage{hyperref}
\hypersetup{colorlinks=true,citecolor=Purple}

\usepackage[capitalise]{cleveref}
\crefformat{equation}{#2(#1)#3}
\theoremstyle{plain}
\newkeytheorem{theorem}[
  parent=section,
  name=Theorem,
  ]
\newkeytheorem{lemma}[
  name=Lemma,
  sibling=theorem,
  Refname={Lemma,Lemmas},
  ]
\newkeytheorem{prop}[
  name=Proposition,
  sibling=theorem,
   Refname={Proposition,Propositions},
  ]
\newkeytheorem{remarks}[
  name=Remarks,
  sibling=theorem,
  style=remark,
  ]

\newcommand{\chg}[1]{
#1}

\begin{document}
\title{Optimal Play in `Guess Who?'}
\author[Cushing]{\vspace{-10pt}David Cushing}
\address{\vspace{-10pt}Dept.~of Mathematics\\
  Alan Turing Building,
  Manchester,
  M13 9PL, UK}
\email{davidcushing1024@gmail.com}
\author[Gipp]{\vspace{-10pt}Stuart Gipp}
\address{\vspace{-10pt}Merry Hell Comic\\
  Trumpington, Cambridge}
\email{stuartgipp@gmail.com}
\author[Levick]{\vspace{-10pt}Ezra Levick}
\address{\vspace{-10pt}Pump Court Chambers\\
  3 Pump Court, London,
  EC4Y 7AJ}
\email{e.macdonald@pumpcourtchambers.com}
\author{\vspace{-10pt}Em Rickinson}
\address{\vspace{-10pt}Speira\\ Sunderland, Tyne \& Wear}
\email{em.rickinson@gmail.com}
\author[Stewart]{\vspace{-10pt}David I.~Stewart}
\address{\vspace{-10pt}Dept.~of Mathematics\\
  Alan Turing Building,
  Manchester,
  M13 9PL, UK}
\email{david.i.stewart@manchester.ac.uk}
\begin{abstract}
  We prove an optimal strategy for the children's game \textit{Guess Who?} assuming the  official rules are in use and that both players ask `classical' questions with bipartite responses. Applying a technique described in [Rabern, B \& Rabern, L 2008, 'A simple solution to the hardest logic puzzle ever', \textit{Analysis}, vol. 68, no. 2, pp.~105-112.] allows for questions with tripartite responses; we explain this innovation and give an optimal strategy for two players applying it.
\end{abstract}
\maketitle
\section{Introduction}\newlength{\theparskip}
\setlength{\theparskip}{\parskip}
\setlength{\baselineskip}{1.2\baselineskip}
In Hasbro's children's game \textit{Guess Who?}, two players begin with the same set of 24 faces in front of them, each with its unique name. Depending on the year of manufacture of the game, players each either surreptitiously choose one character using a sliding frame, or they pick a card from a corresponding set of 24 shuffled cards. The character so picked is now that player's \textit{mystery person} and their job is to guess the opponent's mystery person correctly. If playing the card version, they may immediately eliminate their own mystery person from suspicion, because they know the opponent has a different character. To eliminate suspects, one flips down those characters' images on one's board. Each player knows how many the other has left, but not their names. From the rules:
\begin{quotation}
  Until you're ready to guess who the mystery person is, ask your opponent one question per turn. Each question must have either a \textit{yes} or a \textit{no} answer. For example, you may ask: \textit{Does your person have white hair?} Your opponent must then answer either \textit{yes} or \textit{no}.
\end{quotation}

Taking a guess counts as one turn for a player and immediately ends the game, either in their favour (if the guess is correct), or in their opponent's favour if they lose. Otherwise, they take it in turns to ask questions as above that eliminate suspects and flip down their faces. Eventually one player has eliminated all but one suspect that will they name \textit{on their next turn} and win. Or, at some earlier point any player might try to guess before they are completely sure: for example, you would guess if your opponent has one suspect left since you know that they will win on their next turn and you may as well try your luck.

\ifboxes\begin{wrapfigure}[12]{l}[0pt]{0.55\textwidth}\vspace{-0.60cm}
  \framebox[0.55\textwidth]{
    \parbox{0.53\textwidth}{\vspace{-4pt}\small\subsection*{Previous work} \chg{Our solution of the official game is new, though variants of the game have been considered previously. In \cite{oneill} the game considered uses `merciful guessing'; this means that a wrong guess only ends the turn, not the game. O'Neill finds that the split-in-half technique is in fact optimal. By contrast, in \cite{Nic16}, a player wins immediately on reducing the search space to $1$ suspect. (Note that here a player always wins when they have two suspects left; indeed for this variant one may as well do away with the guessing mechanic entirely.) An optimal strategy for Nica's game looks more like our tripartite version: no matter how many suspects you have left, you always need to pay attention to how many suspects your opponent has; if they have more, then split-in-half; if they have fewer, then aim for the closest power of two which is less than theirs.}}
}
\end{wrapfigure}\fi

The classical play involves asking questions which \textit{always} have an answer; by the rules this is either \textit{yes} or \textit{no} and one is therefore able to split the suspect space into two parts. We thus call those questions \textit{bipartite}. This paper arose after the third author pointed out that by deploying a question containing an embedded paradox---as described in \cite{raberns}---we can improve our chances of winning significantly by asking \emph{tripartite} questions. The following is an example:\bigskip

\begin{quotation}
\vspace{-0.25cm}  Does your person have blond hair\\
  \phantom{mmm}\textbf{OR}\hfill $(*)$\\
  do they have brown hair \textbf{AND} the answer to this question is no?
\end{quotation}

If we were to ask you this question and your person had blond hair, then you would say \textit{yes} 
because the first line succeeds. If they had grey hair, then both parts fail and you would say \textit{no}; 
but if they had brown hair then you would find yourself, in effect, answering: `Is the answer to 
this question \textit{no}?'. You cannot answer honestly, so we may assume that your head explodes---and that explosion can be  treated as a third response.\footnote{\chg{One referee's instinct---which is not unique---was to view our tripartite questions to be in breach of the rules. Our third author, a barrister, has consulted some case law on cheating and responds as follows: $(*)$ clearly has a \textit{yes} answer; it also has a \textit{no} answer; hence it is a question with a \textit{yes} or \textit{no} answer. Moreover, those are the only answers that can be truthfully given. So it appears to fulfil all the requirements made of it by the rules. Hence, if we were to be referred to anyone for the crime of making children's heads explode then it should surely be to the Police rather than toy-maker Hasbro.}} The second and fourth authors have prepared a playable version of the tripartite game: \href{https://www.merryhellcomic.com/conjecture}{merryhellcomic.com/conjecture}. 

\cref{mainthm} and \cref{secthm} below describe \textit{Nash equilibria} with \textit{pure} optimal strategies for playing \textit{Guess Who?} using bipartite and tripartite questions respectively. Less technically, the strategies we provide admit no improvement on their chances of winning, no matter what the other player might plan to do.

Even our optimal strategy in the bipartite case  shows some surprising features: first, it is not optimal always to ask questions that divide the search space by half as closely as possible; second, there are three sporadic exceptional cases that do not obey the generic formula. (Things get even more interesting in the tripartite case.)

\ifboxes\begin{wrapfigure}[17]{r}[0pt]{0.53\textwidth}
  \vspace{-1.1cm}\framebox[0.53\textwidth]{\parbox{0.51\textwidth}{\setlength{\parskip}{\theparskip}\vspace{-8pt}\small\subsection*{For the classroom or outreach activity} \chg{Many games you can buy in the shops have been solved; \textit{Connect Four} is a good example: see \cite{Allis1988AKA, Allen_2010}. But we cannot think of any whose optimal strategy can be expressed as elegantly as that of \textit{Guess Who?}. Moreover, there is nothing in our proof that uses anything more advanced than induction and a recursion relation based on elementary probability; a determined A Level student could have written this paper. Hence \textit{Guess Who?} could guide an introduction to Game Theory to interested students through a `real-world' game. (Such an introduction could be supported with the 2001 film \textit{A Beautiful Mind}---its protagonist John Nash wrote a famously short proving the existence of Nash equilibria.)

      Let us suggest a motivating activity: ask students what they would do if they and their opponent both had four suspects left---expecting most to say that they would ask a question to split the search space into half. Then get them to calculate the probability $P([2,2])$ of winning with that strategy and the probability $P([1,3])$, where instead they split the search space into $[1,3]$. They should find that $P([2,2])=\frac{1}{2}$ and $P([1,3])=\frac{9}{16}$, which accounts for one of our exceptions. Then show them the full strategy. Perhaps follow with a discussion of paradoxes and the (ab)use of tripartite questions in \textit{Guess Who?}, illustrated with \href{https://www.merryhellcomic.com/conjecture}{our playable version}.}}}
\end{wrapfigure}\fi

Before we state the theorems, let us first simplify matters. For any version of the game, note you can always make a question from arbitrary subsets of your remaining characters and ask whether they are in there. So, for a bipartite example one could ask:
\ifboxes\begin{quotation}\hspace{-1cm}\parbox{0.4\textwidth}{\setlength{\baselineskip}{1.2\baselineskip}
  \it Setting $X$ to be the set of people $\{$Ajay, Brenda, Chauntcey, Didi$\}$, is the person in set $X$?}
\end{quotation}\fi
With no prior information, then, the characters are all interchangeable, hence the opponent's mystery person should be assumed to be chosen uniformly at random. In particular, the only relevant information that can influence a player's move at any point is how many suspects \textit{they} have left and how many suspects \textit{their opponent} has left. Let the former be $n$ and the latter $m$; then we say that \textit{the player plays from board state $(n,m)$. }

\section{The official rules with classical questions}
Assume in this section that we only ask bipartite questions.
We play from board state $(n,m)$. A player's decision amounts to an integer $1\leq k\leq \frac{n}{2}$ giving the size $|X|$ of the set $X$ as in the question above---or the integer $0$ that represents using the turn to take a guess. For brevity, let us assume we play against \textit{Joe}. Let $P(n,m)$ be our probability of winning from board state $(n,m)$ when both we and Joe deploy optimal strategies.

Let us treat some easy small cases. Note that the guess decision $0$ is optimal in the two  cases $n=1$ and $m=1$; for if we have $1$ person left then we will guess and win, whereas if Joe has $1$ person left then they will win on their next turn, so we had better guess randomly from our remaining $n$---this gives therefore $P(1,m)=1$ and $P(n,1)=\frac{1}{n}$. If $n=2$ and $m\neq 1$, then a guess is only optimal in case $m=2$, while the decision $1$ is always optimal; for if we have two people left, then a question will leave Joe at state $(m,1)$ with certainty, from which they have only a $\frac{1}{m}$ chance of winning, i.e.~giving \textit{us} a $P(2,m)=\frac{m-1}{m}\geq \frac{1}{2}$ chance of winning. Similarly, one may assume one guesses when $n=3$ if and only if $m=1$; hence $1$ is also an optimal decision at state $(3,m)$ for $m\geq 2$. We have thus verified the following theorem whenever $n\leq 3$.

\begin{theorem}
  Define $$
    \mathcal{S}(n,m):=\begin{cases}
      0,  & \text{if $n=1$ or $m = 1$}  \\

      1 , & \text{if $(n,m) = (4,4)$ }  \\

      3 , & \text{if $(n,m) = (6,4)$ }  \\

      5 , & \text{if $(n,m) = (10,4)$ } \\

      \un & \text{otherwise,}
    \end{cases}$$
  where $\un$ is defined to be $1$ if $n=2$ and is $\left\lfloor\frac{n}{4}\right\rfloor+\left\lfloor\frac{n+1}{4}\right\rfloor$ for $n\geq 3$. Then $\S(n,m)$ is an optimal decision at board state $(n,m)$.\label{mainthm}\end{theorem}

\begin{remarks}
(i) Notice that there are no upper bounds on $n$ and $m$. This means we have in fact solved the game of \textit{Guess Who?} that starts with an arbitrarily large number of suspects. Of course, in the official version, you would start with 23 or 24.

(ii) Define also the integer $$\on:=\left\lfloor\frac{n+2}{4}\right\rfloor+\left\lfloor\frac{n+3}{4}\right\rfloor,$$ and notice that $\un+\on=n$.

(iii) For $n\geq 3$, the first few values of $n_<$ are $(1,2,2,2,3,4,4,4).$ The formula for $n_<$ confirms it is the sequence \href{https://oeis.org/A004524}{A004524} from the OEIS \cite{oeis}---ignoring our modification when $n=2$. 
\end{remarks}

We will proceed with the proof of \cref{mainthm} presently. First we need a couple of simple lemmas.

\begin{lemma}
  \label{lembull}Let $k\geq 2$, so $n\geq 4$. Then
  \cref{nminktab} gives expressions for $\ubull{n-k}$ and $\obull{n-k}$ in terms of $\un,\, \on,\, \uk$ and $\ok$ depending on the values of $n$ and $k$ mod $4$.

  \begin{table}\small
    \begin{center}\renewcommand*{\arraystretch}{1.5}
      \begin{tabular}[t]{|c | c | c | c|}\hline
        $n\mod 4$ & \parbox{1.15cm}{\vspace{5pt}\centering $k\mod 4$,\\$k\geq 3$\vspace{5pt}} & $\ubull{n-k}$ & $\obull{n-k}$ \\\hline
        $0$       & any       & $\un-\ok$     & $\on-\uk$     \\
        $1$       & $0,1$     & $\un-\uk$     & $\on-\ok$     \\
        $1$       & $2,3$     & $\un-\uk-1$   & $\on-\ok+1$   \\
        $2$       & any       & $\un-\uk$     & $\on-\ok$     \\
        $3$       & $0,1$     & $\un-\ok$     & $\on-\uk$     \\
        $3$       & $2,3$     & $\un-\ok+1$   & $\on-\uk-1$\\\hline
      \end{tabular}\qquad\quad
      \begin{minipage}[t]{0.45\textwidth}
        \begin{tabular}[t]{|c | c | c | c|}\hline
        $n$                    & $k$ & $\ubull{n-k}$ & $\obull{n-k}$ \\\hline
        $4$                    & $2$ & $n_<-k_>$     & $n_>-k_<$     \\
        $0,2\,\mathrm{mod}\,4$ & $2$ & $n_<-k_<-1$   & $n_>-k_>+1$   \\
        $1,3\,\mathrm{mod}\,4$ & $2$ & $n_<-k_>$     & $n_>-k_<$.\\\hline
      \end{tabular}
      
\vspace{0.8cm}
        \hspace{-0.25cm}\parbox{1.2\textwidth}{\caption{Tables for $(n-k)_<$ and $(n-k)_>$.\label{nminktab}}}
      \end{minipage}
    \end{center}
  \end{table}
\end{lemma}
\begin{proof}
  Let $k=4k'+k_0$ and $n=4n'+n_0$. Then $\ubull{n-k}=2n'-2k'+\ubull{n_0-k_0}$, and  $\obull{n-k}=2n'-2k'+\obull{n_0-k_0}$. Thus each entry of the table is a simple case-by-case check, left to the interested reader.
\end{proof}

In the proof we use the following notation. From board state $(n,m)$ we make an \textit{optimal} decision $k\leq \frac{n}{2}$, that results in Joe playing from board state $(m,k)$ with probability $\frac{k}{n}$ and from board state $(m,n-k)$ with probability $\frac{n-k}{n}$. We encode all this by saying that Joe is playing from board state $(m,[k:n-k])$, and denote Joe's probability of winning from this point---given that we both subsequently make optimal decisions---as $$P(m,[k:n-k]):=\frac{k}{n}P(m,k)+\frac{n-k}{n}P(m,n-k).$$ Our probability of winning is therefore $1-P(m,[k:n-k])$, and we maximise this value if and only if we minimise $P(m,[k:n-k])$, hence we get \begin{equation}
  P(n,m)=1-P(m,[k:n-k]),\label{fund}
\end{equation} which we call the \textit{fundamental recurrence}---for, if one has already calculated all values of $P(m,k)$ for $k\leq n/2$ then an optimal decision can be read off from those values as any $k$ that minimises $P(m,[n:k])$.

The theorem claims that most of our decisions do not require consideration of the number of suspects left on Joe's board, nor Joe's consideration of our progress. If Joe makes decision $\ell$ indepedently of our $[k:n-k]$ suspects, then they leave us in board state $([k:n-k],[\ell:m-l])$---for example with a $\frac{\ell}{m}\cdot\frac{n-k}{n}$ chance of being in board state $(n-k,\ell)$. For brevity, write $\wtm:=[\ell:m-\ell]$  and note lastly that if we then make decision $[k_1:k_2]$
from state $([k:n-k],\wtm)$, then Joe plays from $(\wtm,[k_1:k-k_1:k_2:n-k-k_2])$ and wins with probability $P(\widetilde{m},[k_1:k-k_1:k_2:n-k-k_2])$.

If our decisions lead to board state $(\wtm,[a_1:\dots:a_r])$ with $\sum a_i=n$, then we are in board state $(\wtm,a_i)$ with probability $\frac{a_i}{n}$. Hence the order of the $a_i$ does not matter and in particular:
\begin{lemma}
  We have $P(\widetilde{m},[a:b:c:d])=P(\widetilde{m},[a:c:b:d])=P(\widetilde{m},[b:d:a:c])$.\label{swap}
\end{lemma}
Lastly, while a proof of \cref{mainthm} entirely by hand would be possible, we have decided that the case-by-case analysis would obscure the main part of the argument. Thus we allow ourselves the luxury of using some Python code\footnote{The code can be found at \url{github.com/cushydom88}} that uses the fundamental recurrence to verify:
\begin{lemma}\label{complem} \cref{mainthm} holds for $n\leq 11$ and $m\leq 4$.\end{lemma}

\begin{proof}[Proof of \cref{mainthm}]
  From \cref{complem} and the remarks before the theorem, we may assume $2\leq m\leq 4$ and  $n\geq 11$, or $m\geq 5$ and $4\leq n\leq 10$. We treat the more tedious cases where $m\leq 4$ and $n\geq 11$ in \cref{lemknot1} after the theorem. So we may assume $m\geq 5$ and $4\leq n\leq 10$.

  It is helpful to prove the the two following statements by induction together:
  \begin{align}
    \begin{split}
      \text{For }m\geq 2,\quad P(m,[k:n-k])
       & \geq P(m,[n_<:n_>])=1-P(n,m)                       \\
       & \text{ whenever }(n,m)\neq (4,4),\ (6,4),\ (10,4),
      \quad\text{and}\end{split}\label{eq1} \\
    \begin{split}
      P([k:n-k],m) & \leq P([n_<:n_>],m)
      \quad\text{ whenever }m\geq 2.
      \label{eq2}
    \end{split}
  \end{align}
  Note that the first line asserts for $k=n_<$ the minimality of $P(m,[k:n-k])$, which verifies the fundamental recurrence, thence the theorem. Assume then that $n+m=d$ and these statements are true whenever $n+m<d$.

  \cref{lemknot1} shows we may assume $k\neq 1$. With that in hand, we start by proving the first line---using the second. Assume $k$ minimises $P(m,[k:n-k])$. Since $m\geq 5$, then by induction using the first line at $(m,n_<)$ and $(m,n_>)$ we get $P(m,[n_<:n_>])=1-P([n_<:n_>],\wtm)$.  And using the second inequality with $m$ replaced by $m_<$ and $m_>$ we get:
  \begin{align}
    \begin{split}1-P([k:n-k],\wtm)\geq 1-P([n_<:n_>],\wtm) & =P(m,[n_<:n_>])                 \\
                                                       & \quad\geq {P(m,[k:n-k])}        \\
                                                       & \qquad \geq {1-P([k:n-k],\wtm)}
    \end{split}\end{align}
  and we have equality throughout.

  Now we must prove the second line. 
  By \cref{lemknot1}(iii) we may assume 
  $(k,m)\neq (4,4),(6,4)$, $(10,4)$. Suppose $P([k:n-k],m)\geq P([n_<:n_>],m)$. Induction yields that $[k_<:(n-k)_<]$ is an optimal decision at $P([k:n-k],m)$, and so
  \begin{align}
    P([n_<:n_>],m)\leq P([k:n-k],m)=1-P(m,[k_<:k_>:(n-k)_<:(n-k)>])
    \label{eq3}
  \end{align}

  Now we need to divide into cases according to \cref{lembull}.

  Suppose first that $(n-k)_<=n_<-k_>$; by \cref{lembull} this happens when $k>2$ and $n=0\mod 4$;  $(n,k)=(3,0)\mod 4$; $(n,k)=(3,1) \mod 4$; or that $k=2$ and either $n=1,3\mod 4$ or $n=4$. We may rewrite the last term of \cref{eq3} to become
  \[P(m,[\uk:\ok:\un-\ok:\on-\uk])= P(m,[\ok:\un-\ok:\uk:\on-\uk])\]
  using \cref{swap}. The expression on the right-hand side is the probability arising from taking decision $[k_>:k_<]$ at board state $([n_<:n_>],m)$, hence is at most $P([n_<:n_>],m)$. It follows that we must have equality all the way through \cref{eq3}, showing that $\S(n,m)=n_<$ is optimal.

  The argument in cases $n=2\mod 4$ or $(n,k)=(3,0),(3,1)\mod 4$ is similar. This time \cref{lembull} gives \begin{align*}P(m,[\uk:(n-k)_<:\ok:(n-k)_>]) & =P(m,[\uk:\un-\uk:\ok:\on-\ok]) \\
                                             & \quad\geq 1-P([\un:\on],m).\end{align*}

  Now assume $(n,k)=(1,2),(1,3)\mod 4$ or that $k=2$ and $n=0,2\mod 4$. Then we have instead from \cref{lembull} that \begin{align*}P(m,[\un:\on])&\geq P(m,[\uk:\un-\uk-1:\ok:\on-\ok+1])\\&\quad \geq 1-P([\un-1:\on+1],m).\end{align*}
  If $n=8n'+1$, then induction yields optimal moves at $(\un-1,m)$ and $(\on+1,m)$ to be $\ubull{\un-1}=2n'-1$ and $\ubull{\on+1}=2n'$, respectively, while those at $(\un,m)$ and $(\on,m)$ are both $2n'$. Thus \begin{align*}P(m,[2n':2n':2n':2n'+1])&\geq P(m,[\un-1:\on+1]\\&\quad =P(m,[2n'-1:2n':2n':2n'+2])\end{align*} and it follows that $P(m,[2n':2n'+1])\geq P(m,[2n'-1:2n'+1])$. However, we may now invoke condition \cref{eq1} with $k=2n'-1$ to get $P(m,[2n'-1:2n'+2])\geq P(m,[2n',2n'+1])$ and hence we must have equality throughout. Thus $\S(n,m)=\un$ was optimal.

  Similarly if $n=8n'+5$, we find $[n_{<<},n_{<>}:n_{><}:n_{>>}]=[2n':2n'+2:2n'+1:2n'+2]$, while $[\ubull{\un-1}:\obull{\un-1}:\ubull{\on+1}:\obull{\on+1}])=[2n':2n'+1:2n'+2:2n'+2]$ and equality of $P(m,[\un:\on])$ and $P(m,[k:n-k])$ is immediate.

  The argument is  even easier when $(n,k)=(3,0),(3,1)\mod 4$. In that case $n_<+1=n_>$ so that
  \[[(n_<+1)_<:(n_<+1)_>:(n_>-1)_<:(n_>-1)_>]=[n_{><}:n_{>>}:n_{<<}:n_{<>}],\] and we are done.
\end{proof}

\newcommand{\vertleq}{\rotatebox{90}{$\,\leq$}}
\newcommand{\vertgeq}{\rotatebox{90}{$\,\geq$}}
\newcommand{\verteq}{\rotatebox{90}{$\,=$}}

\begin{lemma} Suppose $n\geq 4$ and $m\geq 2$. Then:

  (i)  If $(n,m)\neq (4,4)$, then there is an optimal move $k$ at board state $(n,m)$ such that $k>1$.

  (ii) \cref{mainthm} holds for $m\leq 4$. \label{lemknot1}.

  (iii) Statement \cref{eq2} holds for $(k,m)=(4,4),\ (6,4)$ or $(10,4)$.\end{lemma}

\begin{proof}
  Recall the remarks before the theorem that verify \cref{mainthm} for board states $(n,m)$ when $n\leq 3$. \cref{complem} verifies the $n\leq 10$ cases for $m\leq 4$. Thus we assume $n\geq 11$. If $m=2$, then $$P(2,[1:n-1])=\frac{1}{n}\cdot \frac{1}{2}+\frac{n-1}{n}\cdot (1-P(n-1,1))=\frac{1}{2n}+\frac{n-2}{n}=\frac{2n-3}{2n}$$ while for $k\geq 2$, $$P(2,[k:n-k])=1-P([k:n-k],1)=\frac{n-2}{n}<P(2,[1:n-1]).$$
  And if $m=3$, then $$P(3,[1:n-1])=\frac{1}{3n}+\frac{n-1}{n}(1-P(n-1,[1:2]))=\frac{3n-7}{n}$$ whereas for $k\geq 2$, $$P(3,[k:n-k])=1-P([k:n-k],[1:2])=\frac{3n-10}{3n}$$ which is strictly better than $P(3,[1:n-1])$. Lastly suppose $m=4$. Then if $k=1$, we have $$P(4,[1:n-1])=\frac{1}{4n}+\frac{n-1}{n}(1-P(n-1,2))=\frac{4n-11}{4n}.$$ Suppose $k=4$. Then from the computer output in \cref{complem} or by hand, $P(4,4)=\frac{9}{16}$. Hence, $$P(4,[4:n-4])=\frac{9}{4n}+\frac{n-4}{n}(1-P(n-4,2))=\frac{4n-15}{4n}.$$ Now if $k\geq 2$, $k\neq 4$, we have $P(4,[k:n-k])=1-P([k:n-k],2)=\frac{n-4}{n}$ which is strictly smaller than both $P(4,[1:n-1]$ and $P(4,[4:n-4])$. Thus $k\neq 1,4$, and $k=n_<\geq 5$ is optimal, finishing the proof of (i). In fact, we have shown \begin{equation}
    P(n,4)=1-\frac{n-4}{n}=\frac{4}{n}\quad(\text{for }n\geq 10)\label{pn4},
  \end{equation}
  where the case $n=10$ is taken from \cref{complem}.

  We may now assume $m\geq 5$ and we must show $k\neq 1$. Suppose instead $k=1$ is optimal. We have:
  \begin{center}\begin{tikzcd}
      P(m,[2:n-2]) \arrow[r,phantom, "\geq",description ]
      \arrow[d,phantom, "\verteq",description]
      & {1-P(n,m)} \arrow[r,phantom, "=",description ]
      &{P(m,[1:n-1])}\arrow[d,phantom,"\verteq",description]\\
      {\frac{4}{mn}+\frac{n-2}{n}P(m,n-2)} \arrow[r,phantom, "\leq",description ] &
      {\frac{4}{mn}+\frac{n-2}{n}P(m,n-1)} &{\frac{1}{mn}+\frac{n-1}{n}P(m,n-1)}
    \end{tikzcd}\end{center}
  Hence $3/m \geq P(m,n-1)$. The latter is at least $P(m,4)=4/m$ by \cref{pn4}, a contradiction. This finishes the proof of (i) and (ii).

  Finally, we prove (iii). \cref{complem} proves the statement for $k,n-k\leq 10$, and from the computer output or otherwise, $$P(4,4)=\frac{9}{16},\quad P(6,4)=\frac{1}{2}\quad\text{ and }\ P(10,4)=\frac{2}{5}.$$ Assume the statement fails. Then $n_< >k\geq 4$ and $P([k:n-k],4)>P([n_<:n_>],4)$. The former yields $n_<>5$ so $n\geq 10$, thus we may deploy \cref{pn4} to get $$P([4:n-4],4)=\frac{25}{4n},\quad P([6:n-6],4)=\frac{7}{n}\quad\text{and}\quad P([10:n-10],4)=\frac{8}{n}.$$ On the other hand, also from \cref{pn4} we have $P([k:n-k],4)=\frac{8}{n}$ for $k\geq 8$---which we know also happens to hold when $k=10$. This contradiction establishes \cref{eq2} for $k=4,6,10$ and we are done.
\end{proof}

\section{With tripartite questions}
We state an optimal strategy in the case where one can ask tripartite questions with embedded paradoxes of the form:
\begin{quotation}
  Is your person a member of subset $X$\\
  \phantom{mmm}\textbf{OR}\hfill $(\dagger)$\\
  are they a member of subset $Y$ \textbf{AND} the answer to this question is no?
\end{quotation}
\chg{We assume we play from board state $(n,m)$. Note that there is nothing to be gained by asking questions of the form $(\dagger)$ where we allow a person $p$ in the intersection $X\cap Y$: for if the mystery person is indeed $p$, then we $(\dagger)$ would generate a \textit{yes} answer; we can still eliminate the complement of $X$, but this is as good as asking a question with $p$ removed from $Y$; hence we assume that $X$ and $Y$ are disjoint. Up to relabelling characters, then, the only relevant information contained in our decision is a triple of integers $(i,j,k)$ such that $0\leq i\leq \frac{n}{3} \leq j\leq \frac{2n}{3}\leq k$ with $i+j+k=n$, and such that $|X|=i$ and $|Y|=j$; or the symbol $\mathcal{G}$ that represents taking a guess. (The integer $i=0$ occurs only in case $n=2$.)}
\begin{theorem}
  \label{secthm}
  Define for $n\leq 9$,
  $$
    \mathcal{S}(n,m)=\begin{cases}
      \mathcal{G},        & n=1\text{ or }m = 1     \\

      \left(0,1,1\right), & n = 2                   \\

      \left(1,1,1\right), & n = 3                   \\

      \left(1,1,2\right), & n = 4                   \\

      \left(1,2,2\right), & n = 5,\ m \leq 4        \\

      \left(1,1,3\right), & n = 5,\ m\geq 5         \\

      \left(2,2,2\right), & n = 6,\ m \leq 4        \\

      \left(1,1,4\right), & n = 6,\ 5 \leq m \leq 7 \\

      \left(1,2,3\right), & n = 6,\ m \geq 8\end{cases}\quad
    \S(n,m)=\begin{cases}

      \left(2,2,3\right), & n = 7, m \leq 4        \\

      \left(1,3,3\right), & n = 7, m\geq 5         \\

      \left(2,3,3\right), & n = 8, m \leq 4        \\

      \left(1,3,4\right), & n = 8, 5 \leq m \leq 7 \\

      \left(2,3,3\right), & n = 8, m \geq 8        \\

      \left(3,3,3\right), & n = 9, m \leq 4        \\

      \left(1,4,4\right), & n = 9, 5 \leq m \leq 7 \\

      \left(3,3,3\right), & n = 9, m \geq 8.
    \end{cases}$$
  Set $k=\lfloor\log_3 n\rfloor$ and define
  $$n_<:=\sum_{i=0}^{3^k-1}\ff{n+i}{3^{k-1}}\quad, n_\circ:=\sum_{i=3^k}^{2\cdot 3^{k}-1}\ff{n+i}{3^{k-1}},\quad n_>:=\sum_{i=2\cdot 3^{k}}^{3^{k+1}-1}\ff{n+i}{3^{k-1}}.$$ For $n\geq 10$, $m\geq 2$, set
  $$\S(n,m)=\begin{cases}
      ( n_<, n_\circ, n_>)                                & m < 3^{k} + 3^{k-1}    ;                                                                 \\
      \left( 3^{k-1}, 3^{k-1}, n - 2\cdot 3^{k-1} \right) & m \geq 3^{k} + 3^{k-1} , n \leq 3^{k} + 2\cdot 3^{k-1}   ;                               \\
      \left( 3^{k-1}, 3^{k-1} + 2, 3^{k} - 1 \right)      & n = 3^{k} + 2\cdot 3^{k} + 1   ;                                                         \\
      \left( 3^{k-1}, n - 3^{k} -  3^{k-1}, 3^{k} \right) & m \geq 3^{k} + 3^{k-1} ,  3^{k} + 2\cdot 3^{k-1} + 1 < n \leq 3^{k-1} + 2\cdot 3^{k}   ; \\
      \left( 3^{k-1} + 2, 3^{k}-1, 3^{k}  \right)         & n = 3^{k-1} + 2\cdot 3^{k} + 1   ;                                                       \\
      \left(n - 2\cdot 3^{k} , 3^{k}, 3^{k} \right)       & m \geq 3^{k} + 3^{k-1} . n > 3^{k-1} + 2\cdot 3^{k} + 1
    \end{cases}.
  $$
  Then $\S(n,m)$ describes an optimal move at board state $(n,m)$.
\end{theorem}
Unsurprisingly, given the number of cases here, the proof of \cref{secthm} is much more intricate than the bipartite case. We did not find an analogue of \cref{nminktab}, which was what we used to avoid having to calculate formulas for the probabilities of winning from every possible board state; here we take the latter brute force approach instead. We defer a full proof to the Appendix, but we give the rough idea here.
\begin{proof}[Sketch proof.] 
 Let $P(n,m)$ be again the probability of winning from board state $(n,m)$ and for $i+j+k=n$ as in the theorem; let $P([i:j:k],m)$ be the probability of winning if one asks $(\dagger)$ with $|X|=i$, $|Y|=j$. Then we have the following recurrence relation:
  \begin{equation*}
    \begin{split}
    P([i:j:k],m) :&=\frac{i}{n}(1-P(m,i)) + \frac{j}{n}(1-P(m,j)) + \frac{k}{n}(1-P(m,k))\\
    &=1-\frac{i}{n}P(m,i)-\frac{j}{n}P(m,j)-\frac{k}{n}P(m,k).
    \end{split}
  \end{equation*}
  If we play optimally, then $(i,j,k)$ maximises the quantity in this equation---equivalently it minimises $\frac{i}{n}P(m,i)+\frac{j}{n}P(m,j)+\frac{k}{n}P(m,k).$ But now $P(m,\alpha)$ with $\alpha\in\{i,j,k\}$ represents the probability of winning from a board of total size $m+\alpha<m+n$. Hence, we know $P(m,\alpha)$ by induction, and so in principle we can look at all the formulas coming from different choices of $(i,j,k)$ and check that $\S(n,m)$ gives the highest. In practice, generating reasonable formulas so as to compare the effect of the different triples of integers involves fixing $25$ ranges for the pair $(n,m)$ based on the $3$-adic expansions of $n$ and $m$, and checking the inductive step in each case.
\end{proof}

\textbf{Acknowledgements.} We are most grateful to the referees for their careful reading and helpful suggestions for improvement. We would like to thank Fuochi Pizza, Newcastle for hosting some of this research.

\textbf{Declaration of interest.} The authors report that there are no competing interests to declare.

\textbf{Funding.} No funding was received.

\textbf{Author Contribution Statement.} Ezra Levick first suggested the possibility of using tripartite questions in \textit{Guess Who?}. This prompted Davids Cushing and Stewart to discover and prove the optimal strategies given here. Em Rickinson and Stuart Gipp developed the tripartite strategy into the playable game. 

\printbibliography
\appendix
\newpage
\begin{changemargin}{-1cm}{-1cm}
\setlength{\parindent}{0pt}
\setlength{\parskip}{\theparskip}
\setlength{\baselineskip}{1.2\baselineskip}
\thispagestyle{empty}
\markboth{\sc Appendix}{\sc Appendix} 
\section{Proof of the optimal strategy for tripartite \textit{Guess Who?}}

Recall the notation $\S(n,m)$ from \cref{secthm}. We write out  $\mathcal{S}$ for low order cases 
\begin{center}
$\mathbf{10 \leq n \leq 27}$
 
\begin{tabular}{ c|c|c|c|c|c } 
 $n$ & $\mathcal{S}(n,m)$ & $n$ & $\mathcal{S}(n,m)$ & $n$ & $\mathcal{S}(n,m)$   \\ 
  \hline
 10 & (3,3,4) & 16 & (5,5,6) & 22 & (7,7,8)  \\ 
  11 & (3,4,4) & 17 & (5,6,6) & 23 & (7,8,8)  \\ 
 12 & (4,4,4) & 18 & (6,6,6) & 24 & (8,8,8)  \\ 
 13 & (4,4,5) & 19 & (6,6,7) & 25 & (8,8,9)  \\ 
 14 & (4,5,5) & 20 & (6,7,7) & 26 & (8,9,9)  \\ 
 15 & (5,5,5) & 21 & (7,7,7) & 27 & (9,9,9)  \\ 
\end{tabular}
\quad 
\begin{tabular}{ c|c|c|c|c|c } 
 $n$ & $\mathcal{S}(n,m)$ & $n$ & $\mathcal{S}(n,m)$ & $n$ & $\mathcal{S}(n,m)$   \\ 
  \hline
 10 & (3,3,4) & 16 & (3,5,8) & 22 & (5,8,9) \\ 
  11 & (3,3,5) & 17 & (3,5,9) & 23 & (5,9,9)  \\ 
 12 & (3,3,6) & 18 & (3,6,9) & 24 & (6,9,9)  \\ 
 13 & (3,3,7) & 19 & (3,7,9) & 25 & (7,9,9)  \\ 
 14 & (3,3,8) & 20 & (3,8,9) & 26 & (8,9,9)  \\ 
 15 & (3,3,9) & 21 & (3,9,9) & 27 & (9,9,9)  \\ 
\end{tabular}

$m \leq 11$ \hspace{6cm} $m \geq 12$ 
\end{center}

\begin{center}
  $\mathbf{28 \leq n \leq 81.}$ 

\begin{tabular}{ c|c|c|c|c|c } 
 $n$ & $\mathcal{S}(n,m)$ & $n$ & $\mathcal{S}(n,m)$ & $n$ & $\mathcal{S}(n,m)$   \\ 
  \hline
 28 & (9,9,10) & 46 & (15,15,16) & 64 & (21,21,22)  \\ 
  29 & (9,9,11) & 47 & (15,15,17) & 65 & (21,21,23)  \\ 
 30 & (9,9,12) & 48 & (15,15,18) & 66 & (21,21,24)  \\ 
 31 & (9,10,12) & 49 & (15,16,18) & 67 & (21,22,24)  \\ 
 32 & (9,11,12) & 50 & (15,17,18) & 68 & (21,23,24)  \\ 
 33 & (9,12,12) & 51 & (15,18,18) & 69 & (21,24,24)  \\ 
 34 & (10,12,12) & 52 & (16,18,18) & 70 & (22,24,24)  \\ 
  35 & (11,12,12) & 53 & (17,18,18) & 71 & (23,24,24)  \\ 
 36 & (12,12,12) & 54 & (18,18,18) & 72 & (24,24,24)  \\ 
 37 & (12,12,13) & 55 & (18,18,19) & 73 & (24,24,25)  \\ 
 38 & (12,12,14) & 56 & (18,18,20) & 74 & (24,24,26)  \\ 
 39 & (12,12,15) & 57 & (18,18,21) & 75 & (24,24,27)  \\ 
  40 & (12,13,15) & 58 & (18,19,21) & 76 & (24,25,27)  \\ 
  41 & (12,14,15) & 59 & (18,20,21) & 77 & (24,26,27)  \\ 
 42 & (12,15,15) & 60 & (18,21,21) & 78 & (24,27,27)  \\ 
 43 & (13,15,15) & 61 & (19,21,21) & 79 & (25,27,27)  \\ 
 44 & (14,15,15) & 61 & (20,21,21) & 80 & (26,27,27)  \\ 
 45 & (15,15,15) & 63 & (21,21,21) & 81 & (27,27,27)  \\ 
\end{tabular}
\quad
\begin{tabular}{ c|c|c|c|c|c } 
 $n$ & $\mathcal{S}(n,m)$ & $n$ & $\mathcal{S}(n,m)$ & $n$ & $\mathcal{S}(n,m)$   \\ 
  \hline
 28 & (9,9,10) & 46 & (9,11,26) & 64 &  (11,26,27)  \\ 
  29 & (9,9,11) & 47 & (9,11,27) & 65 & (11,27,27)  \\ 
 30 & (9,9,12) & 48 & (9,12,27) & 66 & (12,27,27)  \\ 
 31 & (9,9,13) & 49 & (9,13,27) & 67 & (13,27,27)  \\ 
 32 & (9,9,14) & 50 & (9,14,27) & 68 & (14,27,27)  \\ 
 33 & (9,9,15) & 51 & (9,15,27) & 69 & (15,27,27)  \\ 
 34 & (9,9,16) & 52 & (9,16,27) & 70 & (16,27,27)  \\ 
  35 & (9,9,17) & 53 & (9,17,27) & 71 & (17,27,27)  \\ 
 36 & (9,9,18) & 54 & (9,18,27) & 72 & (18,27,27)  \\ 
 37 & (9,9,19) & 55 & (9,19,27) & 73 & (19,27,27)  \\ 
 38 & (9,9,20) & 56 & (9,20,27) & 74 & (20,27,27)  \\ 
 39 & (9,9,21) & 57 & (9,21,27) & 75 & (21,27,27)  \\ 
  40 & (9,9,22) & 58 & (9,22,27) & 76 & (22,27,27)  \\ 
  41 & (9,9,23) & 59 & (9,23,27) & 77 & (23,27,27)  \\ 
 42 & (9,9,24) & 60 & (9,24,27) & 78 & (24,27,27)  \\ 
 43 & (9,9,25) & 61 & (9,25,27) & 79 & (25,27,27)  \\ 
 44 & (9,9,26) & 61 & (9,26,27) & 80 & (26,27,27)  \\ 
 45 & (9,9,27) & 63 & (9,27,27) & 81 & (27,27,27)  \\ 
\end{tabular}

$m \leq 35$ \hspace{6cm} $m \geq 36$
\end{center}

Let $P(n.m)$ be the probability of winning from board state $(n,m)$ when both players use an optimal strategy, and for $(a,b,c)$ such that $a+b+c = n$, let $P(n,m|(a,b,c))$ denote the probability of winning from board state $(n,m)$ if one asks a question that partitions one's $n$ suspects into sets of sizes $a,b$ and $c$. Then we have 

$$ P(n,m|(a,b,c) ) := \frac{a}{n}(1-P(m,a)) + \frac{b}{n}(1-P(m,b)) + \frac{c}{n}(1-P(m,c)). $$

and so 

$$P(n,m) = \max_{n=a+b+c} P(n,m|(a,b,c)).$$

For convenience, set also $A(n,m) := nm\cdot P(n,m)$ for $n\geq 1$ and $A(0,m) = 0$; set $$S(a,b,c,m) := nm - A(m,a) - A(m,b) - A(m,c).$$

\begin{prop}
Let $n\geq 2$ and $m\geq 2$. Then
$$A(n,m) = \max_{n=a+b+c} S(a,b,c,m)$$
\end{prop}

\begin{proof}

\begin{align*}
A(n,m) & = nm\cdot P(n,m)
\\
& = nm\cdot \max_{n=a+b+c} P(n,m|(a,b,c))
\\
& = nm\cdot \max_{n=a+b+c} \left( \frac{a}{n}(1-P(m,a)) + \frac{b}{n}(1-P(m,b)) + \frac{c}{n}(1-P(m,c))\right)
\\
& = \max_{n=a+b+c} (am(1-P(m,a)) + bm(1-P(m,b)) + cm(1-P(m,c)))
\\
& = \max_{n=a+b+c} (nm - A(m,a) - A(m,b) - A(m,c))
\\
& =  \max_{n=a+b+c} S(a,b,c,m).
\end{align*}
    
\end{proof}

\begin{lemma}\label{nsmol}

The following hold

$$A(n,m) = \begin{cases}
m, & \text{if $n=1$ } \\

1, & \text{if $n>1,$ $m=1$ } \\

2m - 2, & \text{if $n=2,$ $2 \leq  m$ } \\

3m - 3, & \text{if $n=3,$ $2 \leq  m$ } 

		 \end{cases}
$$

\end{lemma}

\begin{proof}

\begin{enumerate}

\item $n=1.$

Guessing guarantees victory. Thus

$$A(1, m) = m\cdot P(1,m) = m.$$

\item $n>1,$ $m=1.$

Asking a question always loses in this case. Therefore the optimal play is making a guess and so

$$A(n,1) =  n\cdot P(n,1) = n\frac{1}{n} = 1.$$

\item $n=2,$ $2 \leq  m.$ 

Guessing at random wins with probability $\frac{1}{2}$  giving $A(2,m) \geq 2m\frac{1}{2} = m.$ 

The only other move is $(1,1).$ In that case

$$S(0,1,1,m) = 2m -  2\cdot A(m, 1) = 2m - 2. $$

Since $2m - 2\geq m$ we are done.

\item $n=3,$ $2 \leq  m.$

Guessing at random wins with probability $\frac{1}{3}$  giving $A(3,m) \geq 3m\frac{1}{3} = m.$ 

It is clear that the move $(1,1,1)$ is at least as good as the move    $(1,2).$ Then

$$S(1,1,1,m) = 3m -  3\cdot A(m, 1) = 3m - 3. $$

Since $3m - 3\geq m$ we are done.\qedhere
\end{enumerate}
\end{proof}

Suppose $1 + 3^{k-1} \leq n \leq 3^{k},$ define
$$\Phi(n) = 3^{k-1} \left(1 + 5n - 3^k - 2 \left\lfloor\frac{n+1}{2}\right\rfloor \right)
$$

Let $k\geq3.$ For $3^{k-1}< n \leq 3^k, 1\leq m$ we define the quantity $B(n,m)$ as follows

\begin{enumerate}

\item $1 + 3^{k-1} \leq n\leq 3^{k},$  $m \leq 3^{k-2}.$ 

$$ B(n, m) = \Phi(m).$$

\item $1 + 3^{k-1} \leq n\leq 3^{k},$  $m = 1 + 3^{k-2}.$

$$ B(n, 1+3^{k-2}) = \begin{cases}
n\cdot 3^{k-2} + 3^{k-2}( 1 +  3^{k-3} )  , & \text{if $1 + 3^{k-1} \leq n \leq 3^{k-1} + 3^{k-2}$ } \\

(1+2\cdot 3^{k-4})n - 3^{k-3}( 1 - 13\cdot 3^{k-4} )   , & \text{if $ 3^{k-1} + 3^{k-2} \leq n \leq  2\cdot 3^{k-1}$ } \\

(-1 + 2\cdot 3^{k-4})n + 3^{k-3}( 11 + 13\cdot 3^{k-4} )  , & \text{if $  2\cdot 3^{k-1} \leq n \leq 2\cdot 3^{k-1}+ 3^{k-2}$ } \\

n - 3^{k-2}( 1 - 3^{k-2} )  , & \text{if $  2\cdot 3^{k-1} + 3^{k-2} \leq n \leq 3^{k}$ } 
\end{cases}
$$

\item $n = 1 + 3^{k-1}.$

$$ B(1+3^{k-1}, m) = \begin{cases}

(1 + 2\cdot 3^{k-2})m - 3^{k-2}( 1 +  3^{k-2} ), & \text{if $ 3^{k-2}  <  m \leq 3^{k-2} + 3^{k-3} $ } \\

  7m\cdot 3^{k-3} + 3^{k-3}( 1 - 13\cdot 3^{k-3} ), & \text{if $ 3^{k-2} + 3^{k-3} \leq m \leq  2\cdot 3^{k-2} $ } \\

(2 + 7\cdot 3^{k-3})m - 3^{k-3}( 11 + 13\cdot 3^{k-3} ), & \text{if $ 2\cdot 3^{k-2} \leq m \leq 2\cdot 3^{k-2} + 3^{k-3} $ } \\

m\cdot 3^{k-1} +  3^{k-2}( 1 - 3^{k-1} ), & \text{if $  2\cdot 3^{k-2} + 3^{k-3}  \leq m \leq   m \leq 2\cdot 3^{k-1} $ } \\

(1 +  3^{k-1})m - 3^{k-3}( 15 + 3^{k} ), & \text{if $  2\cdot 3^{k-1} \leq m  $ }

\end{cases}
$$

\item $ 1 + 3^{k-1} < n \leq  3^{k-1} + 3^{k-2},$  $ 1 +  3^{k-2} < m \leq 3^{k-2} + 3^{k-3}.$

$$ B(n, m) = (n - 3^{k-2})m - 3^{k-2}(n - 2\cdot 3^{k-2}). $$

\item $ 1 + 3^{k-1} < n \leq  3^{k-1} + 3^{k-2},$  $ 3^{k-2} + 3^{k-3} \leq m \leq 3^{k-2} + 2\cdot 3^{k-3}.$

$$B(n,m) = \left(n-1 + 3^{k-3} - \left\lfloor\frac{n-1}{3}\right\rfloor \right)m + 4\cdot 3^{k-3}\left\lfloor\frac{n-1}{3}\right\rfloor + 3^{k-3}( -3n + 4  + 2\cdot 3^{k-3}).$$

\item $ 1 + 3^{k-1} < n \leq  3^{k-1} + 3^{k-2},$  $ 3^{k-2} + 2\cdot 3^{k-3} \leq m \leq 2\cdot 3^{k-2} .$

$$B(n,m) = \left(n-1 + 3^{k-3} - \left\lfloor\frac{n-1}{3}\right\rfloor \right)m + 4\cdot 3^{k-3}\left\lfloor\frac{n-1}{3}\right\rfloor + 3^{k-3}( -3n + 4  + 2\cdot 3^{k-3}).$$

\item $ 1 + 3^{k-1} < n \leq  3^{k-1} + 3^{k-2},$  $  2\cdot 3^{k-2} \leq m \leq  2\cdot 3^{k-2} + 3^{k-3} .$

$$B(n,m) = \left(2 + 3^{k-3} + 2\left\lfloor\frac{n-1}{3}\right\rfloor \right)m  - 14\cdot 3^{k-3}\left\lfloor\frac{n-1}{3}\right\rfloor + 3^{k-3}( 3n  - 14  + 2\cdot 3^{k-3}).$$

\item $ 1 + 3^{k-1} < n \leq  3^{k-1} + 3^{k-2},$  $   2\cdot 3^{k-2} + 3^{k-3} \leq m \leq 2\cdot 3^{k-2} + 2\cdot 3^{k-3}.$

$$B(n,m) = m\cdot 3^{k-1} + 3^{k-2}( n - 2\cdot 3^{k-1}). $$

\item $ 1 + 3^{k-1} < n \leq  3^{k-1} + 3^{k-2},$  $  2\cdot 3^{k-2} + 2\cdot 3^{k-3} \leq m \leq 3^{k-1}.$

$$B(n,m) = m\cdot 3^{k-1} + 3^{k-2}( n - 2\cdot 3^{k-1}). $$

\item $ 3^{k-1} + 3^{k-2} < n \leq  3^{k-1} + 2\cdot 3^{k-2},$  $1+3^{k-2} <  m \leq 3^{k-1}.$

$$B(n,m) = m\cdot 3^{k-1} +  \left\lfloor\frac{m-1}{3}\right\rfloor(n - 4\cdot 3^{k-2}  ) +  n  - 3^{k-3}\left( n + 12 + 2\cdot 3^{k-2} \right) $$

\item $   3^{k-1} + 2\cdot 3^{k-2} \leq n \leq  2\cdot 3^{k-1},$  $1+ 3^{k-2} <  m \leq 3^{k-1}.$

$$B(n,m) = m\cdot 3^{k-1} +  \left\lfloor\frac{m-1}{3}\right\rfloor(n - 4\cdot 3^{k-2}  ) +  n  - 3^{k-3}\left( n + 12 + 2\cdot 3^{k-2} \right) $$

\item $  2\cdot 3^{k-1} \leq n \leq 2\cdot 3^{k-1} + 3^{k-2},$  $1+3^{k-2} <  m \leq 3^{k-1}.$

$$B(n,m) = (n - 3^{k-1} )m +\left\lfloor\frac{m-1}{3}\right\rfloor(14\cdot 3^{k-2}  -2n) - 2n - 3^{k-3} \left( n - 42 + 2\cdot 3^{k-2}  \right). $$

\item $ 2\cdot 3^{k-1} + 3^{k-2} \leq n \leq 2\cdot 3^{k-1} + 2\cdot 3^{k-2} ,$  $1+ 3^{k-2} <  m \leq 3^{k-1}.$

$$B(n,m) = (n - 3^{k-1} )m  + \left\lfloor\frac{m-1}{2}\right\rfloor( 14\cdot 3^{k-2} - 2n  ) - n + 3^{k-2} \left( 7 - 3^{k-1}\right).  $$

\item $2\cdot 3^{k-1} + 2\cdot 3^{k-2} \leq n \leq 3^{k},$  $1+3^{k-2} <  m \leq 3^{k-1}.$

$$B(n,m) = m\cdot 3^{k-1} + 2\cdot 3^{k-2}\left\lfloor\frac{m}{2}\right\rfloor + 3^{k-2} - 3^{2k-3} . $$

\item $ 1 + 3^{k-1} < n \leq 3^{k-1} + 3^{k-2},$  $1+ 3^{k-1} <  m < 3^{k-1} + 3^{k-2}.$

$$B(n,m) = m\cdot 3^{k-1} + 3^{k-2}(n - 2\cdot 3^{k-1} ).$$

\item $ 3^{k-1} + 3^{k-2} \leq n \leq 3^{k-1} + 2 \cdot 3^{k-2},$  $1+ 3^{k-1} \leq  m < 3^{k-1} + 3^{k-2}.$

$$B(n,m) =m\cdot 3^{k-1}+ \left\lfloor\frac{m-1}{3}\right\rfloor( n -4\cdot 3^{k-2}   )  +  n  - 3^{k-3}\left( n  + 12 + 2\cdot 3^{k-2} \right).$$

\item $ 3^{k-1} + 2 \cdot 3^{k-2} \leq n \leq 2 \cdot 3^{k-1},$  $1+3^{k-1} \leq  m < 3^{k-1} + 3^{k-2}.$

$$B(n,m) =m\cdot 3^{k-1}+ \left\lfloor\frac{m-1}{3}\right\rfloor( n -4\cdot 3^{k-2}   )  +  n  - 3^{k-3}\left( n  + 12 + 2\cdot 3^{k-2} \right).$$

\item $  2 \cdot 3^{k-1} \leq n \leq  2 \cdot 3^{k-1} + 3^{k-2},$  $1+3^{k-1} \leq  m < 3^{k-1} + 3^{k-2}.$

$$B(n,m) = (n -  3^{k-1}  )m  + \left\lfloor\frac{m-1}{3}\right\rfloor (14\cdot 3^{k-2} - 2n )  - 2n + 3^{k-3} \left( - n + 42 - 2\cdot 3^{k-2}  \right).  $$

\item $  2 \cdot 3^{k-1} + 3^{k-2} \leq n \leq 2 \cdot 3^{k-1} + 2\cdot 3^{k-2} ,$  $1 + 3^{k-1} \leq  m < 3^{k-1} + 3^{k-2}.$

$$B(n,m) = (n - 3^{k-1}  )m + 3^{k-1}( 2\cdot 3^{k-1} - n  ). $$

\item $ 2 \cdot 3^{k-1} + 2\cdot 3^{k-2} \leq n \leq 3^{k},$  $1+3^{k-1} \leq  m < 3^{k-1} + 3^{k-2}.$

$$B(n,m) = (n - 3^{k-1}  )m + 3^{k-1}( 2\cdot 3^{k-1} - n  ). $$

\item $ 1 + 3^{k-1} < n \leq 3^{k},$  $3^{k-1} + 3^{k-2} \leq  m \leq 3^{k-1} + 2\cdot 3^{k-2} .$

$$B(n,m) = \left(n + 3^{k-2} - 1 - \left\lfloor\frac{n-1}{3}\right\rfloor \right)m    + 3^{k-2}\left( - 3n + 4 + 2\cdot 3^{k-2} + 4\left\lfloor\frac{n-1}{3}\right\rfloor \right).$$

\item $ 1 + 3^{k-1} < n \leq 3^{k},$  $3^{k-1} + 2\cdot 3^{k-2} \leq m \leq 2\cdot 3^{k-1}  .$

$$B(n,m) = \left(n + 3^{k-2} - 1 -  \left\lfloor\frac{n-1}{3}\right\rfloor \right)m    + 3^{k-2}\left( - 3n + 4 + 2\cdot 3^{k-2} + 4\left\lfloor\frac{n-1}{3}\right\rfloor \right).$$

\item $ 1 + 3^{k-1} < n \leq 3^{k},$  $2\cdot 3^{k-1} \leq m \leq 2\cdot 3^{k-1} + 3^{k-2} .$

$$B(n,m) = \left(2 + 3^{k-2} + 2\left\lfloor\frac{n-1}{3}\right\rfloor  \right)m  + 3^{k-2}\left( 3n - 14 + 2\cdot 3^{k-2} - 14\left\lfloor\frac{n-1}{3}\right\rfloor    \right) $$

\item $ 1 + 3^{k-1} < n \leq 3^{k},$  $ 2\cdot 3^{k-1} + 3^{k-2} \leq m \leq 2\cdot 3^{k-1} + 2\cdot 3^{k-2} .$

$$B(n,m) = \left(1 + 2\left\lfloor\frac{n-1}{2}\right\rfloor\right)m + 3^{k-2} \left( 3n - 7 + 3^{k} - 14\left\lfloor\frac{n-1}{2}\right\rfloor \right). $$

\item $ 1 + 3^{k-1} < n \leq 3^{k},$  $  2\cdot 3^{k-1} + 2\cdot 3^{k-2} \leq m .$

$$B(n,m) = nm + 2\cdot 3^{k-2}\left\lfloor\frac{n-1}{2}\right\rfloor - 3^{k-2}( 5n -1 -  3^{k} ).  $$

\end{enumerate}

\begin{lemma}\label{bound_lemma}
For all $n,m$ we have 
$$A(n,m+1) - A(n,m) \leq n.$$
\end{lemma}

\begin{proof}
Checking by computer we see the result holds for $m,n <= 9.$

Let $(a,b,c)$ be an optimal play for the state $(n,m+1).$ Since $A(n,m) \geq nm - A(m,a) - A(m,b) - A(m,c)$ we have
    $$A(n,m+1) - A(n,m) \leq m + A(m,a) - A(m+1,a) + A(m,b) -A(m+1,b) + A(m,c) - A(m+1,c).$$

Thus it is enough to show $A(m,s) - A(m+1,s) \leq 0 $ for $s= a,b,c.$

Let $(x,y,z)$ be an optimal strategy for the state $(m,a).$ Using $A(m+1,a) \geq ma + a - A(a,x) - A(a,y) - A(a,z+1)$ we have

$$A(m,a) - A(m+1,a) \leq A(a,z+1) - A(a,z) - a.$$

Thus we have reduced the problem to show that $A(a,z+1)-A(a,z)\leq a.$ Since $a < n$ and $z < m$ we repeat this approach until we reach small enough values where the result is verified.
\end{proof}


\begin{lemma}
Let $k \geq 4$ and suppose that $A(n',m') = B(n',m')$  for all $n' \leq 3^{k-1}.$ Let $3^{k-1} < n \leq 3^{k}. $  Then $A(n,m) = B(n,m) = \Phi(m)$ for all $m \leq 3^{k-2}.$
\end{lemma}

\begin{proof}
Suppose $ 3^{\ell - 1} < m \leq 3^{\ell} \leq 3^{k-2}.$

For $x \geq  3^{\ell}$ we have

$$A(m, x) = B(m, x) =  mx + 3^{\ell-2}\left( -5m - 1 +  3^\ell  + 2\left\lfloor\frac{m + 1 }{2}\right\rfloor \right) = mx  - \frac{\Phi(m)}{3}. $$

Suppose $(a,b, n - a - b)$ is an optimal play for the state $(n, m)$ with $a\leq b \leq n - a - b.$ Thus $n - a - b \geq  3^{k-2} \geq  3^{\ell}$ and we have

$$A(n, m) = S(a,b,c,m) = am - A(m,a) + bm - A(m,b) + \frac{\Phi(m)}{3}.$$

By applying \cref{bound_lemma} repeatedly we have
$$am - A(m,a) \leq am + m - A(m,a+1) = (a+1)m  - A(m,a+1) \leq \cdots \leq m \cdot 3^{k} - A(m, 3^k) $$

 $$ =m \cdot 3^{k} -  m\cdot 3^{k} + 3^{\ell-2}\left( 5m + 1 -  3^\ell  - 2\left\lfloor\frac{m + 1 }{2}\right\rfloor \right)  = \frac{\Phi(m)}{3}.$$

Likewise $bm - A(m,b) \leq \frac{\Phi(m)}{3}.$  Thus $A(n,m) \leq \Phi(m).$ 

Choose $a', b', c'$ such that $3^{k-2} \leq a' \leq b' \leq c'$ and $a' + b' + c' = n.$ Then

$$A(n,m) \geq S(a',b',c',m) = nm - a'm + \frac{\Phi(m)}{3} - b'm + \frac{\Phi(m)}{3} - c'm + \frac{\Phi(m)}{3} = \Phi(m),$$

completing the proof.
\end{proof}

Let $k\geq 4$. In the following we will assume that $A(n,m) = B(n,m)$ for all $n \leq 3^{k-1}.$ For $m\leq 3^{k-1}$ we are therefore assuming that the following table of induction formula hold. 

\begin{enumerate}

\item $1 + 3^{k-2} < m $, $x\leq 3^{k-3}$

$$A(m,x) = \Phi(x).$$

\item $m = 1+3^{k-2}.$

$$ A(1+3^{k-2}, x) = \begin{cases}

(1 + 2\cdot 3^{k-3})x - 3^{k-3}( 1 +  3^{k-3} ), & \text{if $ 3^{k-3} \leq  x \leq 3^{k-2} + 3^{k-4} $ } \\

  7x\cdot 3^{k-4} + 3^{k-4}( 1 - 13\cdot 3^{k-4} ), & \text{if $ 3^{k-3} + 3^{k-4} \leq x \leq  2\cdot 3^{k-3} $ } \\

(2 + 7\cdot 3^{k-4})x - 3^{k-4}( 11 + 13\cdot 3^{k-4} ), & \text{if $ 2\cdot 3^{k-3} \leq x \leq 2\cdot 3^{k-3} + 3^{k-4} $ } \\

x\cdot 3^{k-2} +  3^{k-3}( 1 - 3^{k-2} ), & \text{if $  2\cdot 3^{k-3} + 3^{k-3}  \leq z \leq 3^{k-2} $ } 

\end{cases}
$$

\item $ 1 + 3^{k-2} < m \leq  3^{k-2} + 3^{k-3},$  $ 3^{k-3} < x \leq 3^{k-3} + 3^{k-4}.$

$$ A(m,x) = (m - 3^{k-3})x - 3^{k-3}(m - 2\cdot 3^{k-3}). $$

\item $ 1 + 3^{k-2} < m \leq  3^{k-2} + 3^{k-3},$  $ 3^{k-3} + 3^{k-4} \leq x \leq 3^{k-3} + 2\cdot 3^{k-4}.$

$$A(m,x) = \left(m-1 + 3^{k-4} - \left\lfloor\frac{m-1}{3}\right\rfloor \right)x + 4\cdot 3^{k-4}\left\lfloor\frac{m-1}{3}\right\rfloor + 3^{k-4}( -3m + 4  + 2\cdot 3^{k-4}).$$

\item $ 1 + 3^{k-2} < m \leq  3^{k-2} + 3^{k-3},$  $ 3^{k-3} + 2\cdot 3^{k-4} \leq x \leq 2\cdot 3^{k-3} .$

$$A(m,x) = \left(m-1 + 3^{k-4} - \left\lfloor\frac{m-1}{3}\right\rfloor \right)x + 4\cdot 3^{k-4}\left\lfloor\frac{m-1}{3}\right\rfloor + 3^{k-4}( -3m + 4  + 2\cdot 3^{k-4}).$$

\item $ 1 + 3^{k-2} < m \leq  3^{k-2} + 3^{k-3},$  $  2\cdot 3^{k-3} \leq x \leq  2\cdot 3^{k-3} + 3^{k-4} .$

$$A(m,x) = \left(2 + 3^{k-4} + 2\left\lfloor\frac{m-1}{3}\right\rfloor \right)x  - 14\cdot 3^{k-4}\left\lfloor\frac{m-1}{3}\right\rfloor + 3^{k-4}( 3m  - 14  + 2\cdot 3^{k-4}).$$

\item $ 1 + 3^{k-2} < m \leq  3^{k-2} + 3^{k-3},$  $   2\cdot 3^{k-3} + 3^{k-4} \leq x \leq 2\cdot 3^{k-3} + 2\cdot 3^{k-4}.$

$$A(m,x) = x\cdot 3^{k-2} + 3^{k-3}( m - 2\cdot 3^{k-2}). $$

\item $ 1 + 3^{k-2} < m \leq  3^{k-2} + 3^{k-3},$  $  2\cdot 3^{k-3} + 2\cdot 3^{k-4} \leq x \leq 3^{k-2}.$

$$A(m,x) = x\cdot 3^{k-2} + 3^{k-3}( m - 2\cdot 3^{k-2}). $$

\item $ 3^{k-2} + 3^{k-3} < m \leq  3^{k-2} + 2\cdot 3^{k-3},$  $3^{k-3} <  x \leq 3^{k-2}.$

$$A(m,x) = x\cdot 3^{k-2} +  \left\lfloor\frac{x-1}{3}\right\rfloor(m - 4\cdot 3^{k-3}  ) +  m  - 3^{k-4}\left( m + 12  + 2\cdot 3^{k-3} \right)  $$

\item $   3^{k-2} + 2\cdot 3^{k-3} \leq m \leq  2\cdot 3^{k-2},$  $3^{k-3} <  x \leq 3^{k-2}.$

$$A(m,x) = x\cdot 3^{k-2} +  \left\lfloor\frac{x-1}{3}\right\rfloor(m - 4\cdot 3^{k-3}  ) +  m  - 3^{k-4}\left( m + 12 + 2\cdot 3^{k-3} \right).  $$

\item $  2\cdot 3^{k-2} \leq m \leq 2\cdot 3^{k-2} + 3^{k-3},$  $3^{k-3} <  x \leq 3^{k-2}.$

$$A(m,x) = (m - 3^{k-2} )x +\left\lfloor\frac{x-1}{3}\right\rfloor(14\cdot 3^{k-3}  -2m) - 2m - 3^{k-4} \left( m -42 + 2\cdot 3^{k-3}  \right). $$

\item $ 2\cdot 3^{k-2} + 3^{k-3} \leq m \leq 2\cdot 3^{k-2} + 2\cdot 3^{k-3} ,$  $3^{k-3} <  x \leq 3^{k-2}.$

$$A(m,x) = (m - 3^{k-2} )x  + \left\lfloor\frac{x-1}{2}\right\rfloor( 14\cdot 3^{k-3} - 2m  ) - m + 3^{k-3} \left( 7 - 3^{k-2}\right).  $$

\item $2\cdot 3^{k-2} + 2\cdot 3^{k-3} \leq m \leq 3^{k-1},$  $3^{k-3} <  x \leq 3^{k-2}.$

$$A(m,x) = x\cdot 3^{k-2} + 2\cdot 3^{k-3}\left\lfloor\frac{x}{2}\right\rfloor + 3^{k-3} - 3^{2k-5} . $$

\item $ 1 + 3^{k-2} < m \leq 3^{k-2} + 3^{k-3},$  $3^{k-2} <  x < 3^{k-2} + 3^{k-3}.$

$$A(m,x) = x\cdot 3^{k-2} + 3^{k-3}(m - 2\cdot 3^{k-2} ).$$

\item $ 3^{k-2} + 3^{k-3} \leq m \leq 3^{k-
2} + 2 \cdot 3^{k-3},$  $3^{k-2} <  x < 3^{k-2} + 3^{k-3}.$

$$A(m,x) =x\cdot 3^{k-2}+ \left\lfloor\frac{x-1}{3}\right\rfloor( m -4\cdot 3^{k-3}   )  +  m  - 3^{k-4}\left( m  + 12 + 2\cdot 3^{k-3} \right).$$

\item $ 3^{k-2} + 2 \cdot 3^{k-3} \leq m \leq 2 \cdot 3^{k-2},$  $3^{k-2} <  x < 3^{k-2} + 3^{k-3}.$

$$A(m,x) =x\cdot 3^{k-2}+ \left\lfloor\frac{x-1}{3}\right\rfloor( m -4\cdot 3^{k-3}   )  +  m  - 3^{k-4}\left( m  + 12 + 2\cdot 3^{k-3} \right).$$

\item $  2 \cdot 3^{k-2} \leq m \leq  2 \cdot 3^{k-2} + 3^{k-3},$  $3^{k-2} <  x < 3^{k-2} + 3^{k-3}.$

$$A(m,x) = (m -  3^{k-2}  )x  + \left\lfloor\frac{x-1}{3}\right\rfloor (14\cdot 3^{k-3} - 2m )  - 2m + 3^{k-4} \left( - m + 42 - 2\cdot 3^{k-3}  \right).  $$

\item $  2 \cdot 3^{k-2} + 3^{k-3} \leq m \leq 2 \cdot 3^{k-2} + 2\cdot 3^{k-3} ,$  $3^{k-2} <  x < 3^{k-2} + 3^{k-3}.$

$$A(m,x) = (m - 3^{k-2}  )x + 3^{k-2}( 2\cdot 3^{k-2} - m  ). $$

\item $ 2 \cdot 3^{k-2} + 2\cdot 3^{k-3} \leq m \leq 3^{k-1},$  $3^{k-2} <  x < 3^{k-2} + 3^{k-3}.$

$$A(m,x) = (m - 3^{k-2}  )x + 3^{k-2}( 2\cdot 3^{k-2} - m  ). $$

\item $ 1 + 3^{k-2} < m \leq 3^{k-1},$  $3^{k-2} + 3^{k-3} \leq  x \leq 3^{k-2} + 2\cdot 3^{k-3} .$

$$A(m,x) = \left(m + 3^{k-3} - 1 -\left\lfloor\frac{m-1}{3}\right\rfloor \right)x    + 3^{k-3}\left( - 3m + 4 + 2\cdot 3^{k-3} + 4\left\lfloor\frac{m-1}{3}\right\rfloor \right).$$

\item $ 1 + 3^{k-2} < m \leq 3^{k-1},$  $3^{k-2} + 2\cdot 3^{k-3} \leq x \leq 2\cdot 3^{k-2}  .$

$$A(m,x) = \left(m + 3^{k-3} - 1 - \left\lfloor\frac{m-1}{3}\right\rfloor \right)x    + 3^{k-3}\left( - 3m + 4 + 2\cdot 3^{k-3} + 4\left\lfloor\frac{m-1}{3}\right\rfloor \right).$$

\item $ 1 + 3^{k-2} < m \leq 3^{k-1},$  $2\cdot 3^{k-2} \leq x \leq 2\cdot 3^{k-2} + 3^{k-3} .$

$$A(m,x) = \left(2 + 3^{k-3} + 2\left\lfloor\frac{m-1}{3}\right\rfloor  \right)x  + 3^{k-3}\left( 3m - 14 + 2\cdot 3^{k-3} - 14\left\lfloor\frac{m-1}{3}\right\rfloor    \right) $$

\item $ 1 + 3^{k-2} < m \leq 3^{k-1},$  $ 2\cdot 3^{k-2} + 3^{k-3} \leq x \leq 2\cdot 3^{k-2} + 2\cdot 3^{k-3} .$

$$A(m,x) = \left(1 + 2\left\lfloor\frac{m-1}{2}\right\rfloor\right)x + 3^{k-3} \left( 3m - 7 +  3^{k-1} - 14\left\lfloor\frac{m-1}{2}\right\rfloor \right). $$

\item $ 1 + 3^{k-2} < m \leq 3^{k-1},$  $  2\cdot 3^{k-2} + 2\cdot 3^{k-3} \leq x .$

$$A(m,x) = mx + 2\cdot 3^{k-3}\left\lfloor\frac{m-1}{2}\right\rfloor - 3^{k-3}( 5m -1 -  3^{k-1} ).  $$

\end{enumerate}

\begin{lemma}\label{philemma}
    Let $1 + 3^{k-1} < m \leq 3^k.$ We have
    $$\Phi(m)-\Phi(m-1) \leq 5\cdot 3^{k-1}.$$
\end{lemma}

\begin{proof}

\begin{align*}
\Phi(m)-\Phi(m-1) & = 3^{k-1}\left(1 + 5m - 3^k - 2 \left\lfloor\frac{m+1}{2}\right\rfloor \right) - 3^{k-1}\left(1 + 5m - 5 - 3^k - 2 \left\lfloor\frac{m}{2}\right\rfloor \right)
\\
& = 3^{k-1}\left( 5 + 2 \left\lfloor\frac{m}{2}\right\rfloor - 2 \left\lfloor\frac{m+1}{2}\right\rfloor \right)
\\
& \leq 5\cdot 3^{k-1}.
\end{align*}

\end{proof}

\begin{lemma}\label{allgaps}
Set $$G(m,x) = B(m,x+1)- B(m,x).$$

Let $k \geq 4.$ The following hold
\begin{enumerate}
\item Let $ 1 + 3^{k-2} < m \leq 3^{k-1}$ and $3^{k-2} + 3^{k-3}< x \leq y. $ Then
    $$G(m,y-1) \geq G(m,x-1)\geq \frac{2m}{3} + 3^{k-3} - 1 \geq 3^{k-2}.$$

\item Let  $3^{k-2} < m \leq 3^{k-1}$ and $x \leq 3^{k-3}. $ Then
    $$G(m,x) \leq 3^{k-2}.$$

\item Let  $3^{k-2} < m \leq 3^{k-2} + 3^{k-3}$ and $x \leq 3^{k-2}. $ Then
    $$ G(m,x) \leq 3^{k-2}.$$

\item Let  $3^{k-2} < m \leq 3^{k-1}$ and $y > 3^{k-2}. $ Then
    $$ G(m,y-1) \geq 3^{k-2}.$$

\item Let $x\geq 2.$ Then

$$G(2, x)= 2, G(3,x) = 3.$$

\end{enumerate}
    
\end{lemma}

\begin{proof}

\begin{enumerate}
\item  We have 
$$    G(m,c-1) = \begin{cases}

m + 3^{k-3} - 1 -\left\lfloor\frac{m-1}{3}\right\rfloor  , &\text{if $3^{k-2} + 3^{k-3}  < c \leq 2\cdot 3^{k-2}   $ } \\

 2 + 3^{k-3} + 2\left\lfloor\frac{m-1}{3}\right\rfloor, &\text{if $2\cdot 3^{k-2}  < c \leq 2\cdot 3^{k-2} + 3^{k-3} $ } \\

  1 + 2\left\lfloor\frac{m-1}{2}\right\rfloor , &\text{if $ 2\cdot 3^{k-2} +  3^{k-3}  < c \leq 2\cdot 3^{k-2}  + 2\cdot 3^{k-3}  $ } \\

 m, &\text{if $ 2\cdot 3^{k-2} + 2\cdot 3^{k-3} < c$ } 
\end{cases}$$

It suffices to show that
$$2^{k-2}\leq\frac{2m}{3} + 3^{k-3} - 1\leq m + 3^{k-3} - 1 -\left\lfloor\frac{m-1}{3}\right\rfloor\leq  2 + 3^{k-3} + 2\left\lfloor\frac{m-1}{3}\right\rfloor \leq 1 + 2\left\lfloor\frac{m-1}{2}\right\rfloor \leq m.$$

Indeed we have
$$\frac{2m}{3} + 3^{k-3} - 1 \geq \frac{2(3^{k-2}+2)}{3} + 3^{k-3} - 1 \geq 3^{k-2}, $$
$$m + 3^{k-3} - 1 -\left\lfloor\frac{m-1}{3}\right\rfloor \geq m + 3^{k-3} - 1 - \frac{m-1}{3} \geq \frac{2m}{3} + 3^{k-3} - 1,$$

$$\left(  2 + 3^{k-3} + 2\left\lfloor\frac{m-1}{3}\right\rfloor\right) - \left( m + 3^{k-3} - 1 -\left\lfloor\frac{m-1}{3}\right\rfloor \right) = 3\left\lfloor\frac{m-1}{3}\right\rfloor + 3 - m \geq 0,$$

\begin{align*}
 2 + 3^{k-3} + 2\left\lfloor\frac{m-1}{3}\right\rfloor \leq & 2 + 3^{k-3} + 2\frac{m-1}{3}
 \\
\leq & 2 + 3^{k-3} + 2\frac{m-1}{3}
 \\
\leq & \frac{m}{3} - 4 + 2\frac{m-1}{3}
 \\
\leq & m - 2
 \\
\leq &  1 + 2\left\lfloor\frac{m-1}{2}\right\rfloor,
\end{align*}

$$  1 + 2\left\lfloor\frac{m-1}{2}\right\rfloor \leq  1 + 2\frac{m-1}{2} = m - 1 \leq m.$$

\item Let $r$ be such that $ 3^{r-1} <x \leq 3^{r}.$ Note $r\leq k - 3.$ Then by \cref{philemma} we have

    $$ G(m,x) \leq 5\cdot 3^{r-1} \leq 5\cdot 3^{k-4} \leq 3^{k-2}.$$

\item For $m = 1 + 3^{k-2}$ we have

$$ G(1+3^{k-2}, x )= \begin{cases}

1 + 2\cdot 3^{k-3} \leq 3^{k-2}, & \text{if $ 3^{k-3} \leq  x < 3^{k-3} + 3^{k-4} $ } \\

  7\cdot 3^{k-4} \leq 3^{k-2} , & \text{if $ 3^{k-3} + 3^{k-4} \leq x <  2\cdot 3^{k-3} $ } \\

2 + 7\cdot 3^{k-4}\leq 3^{k-2} , & \text{if $ 2\cdot 3^{k-3} \leq x < 2\cdot 3^{k-3} + 3^{k-4} $ } \\

3^{k-2}, & \text{if $  2\cdot 3^{k-3} + 3^{k-2}  \leq x < 3^{k-2} $ } 

\end{cases}
$$

For  $ 1 + 3^{k-2} < m \leq  3^{k-2} + 3^{k-3}$ and $ 3^{k-3} \leq x < 3^{k-2} $ we have

$$ G(m, x )= \begin{cases}

m - 3^{k-3}, & \text{if  $ 3^{k-3} < x \leq 3^{k-3} + 3^{k-4}$ } \\

m-1 + 3^{k-4} - \left\lfloor\frac{m-1}{3}\right\rfloor , & \text{if $ 3^{k-3} + 3^{k-4} \leq x \leq 2\cdot 3^{k-3}$ } \\

2 + 3^{k-4} + 2\left\lfloor\frac{m-1}{3}\right\rfloor , & \text{if $  2\cdot 3^{k-3} \leq x \leq  2\cdot 3^{k-3} + 3^{k-4}$ } \\

3^{k-2}, & \text{if  $   2\cdot 3^{k-3} + 3^{k-4} \leq x \leq 3^{k-2}$ } 

\end{cases}
$$

All of these expressions are all monotonic non-decreasing with respect to $m.$  Thus by using $ m \leq  3^{k-2} + 3^{k-3}$ we have

$$m - 3^{k-3} \leq 3^{k-2} + 3^{k-3} - 3^{k-3} = 3^{k-2},$$

\begin{align*}
& m-1 + 3^{k-4} - \left\lfloor\frac{m-1}{3}\right\rfloor \leq  3^{k-2} + 3^{k-3}-1 + 3^{k-4} - \left\lfloor\frac{3^{k-2} + 3^{k-3}-1}{3}\right\rfloor 
\\
= & 3^{k-2} + 3^{k-3}-1 + 3^{k-4} - (3^{k-2} + 3^{k-3} - 1) =  3^{k-2},
\end{align*}

and

\begin{align*}
& 2 + 3^{k-4} + 2\left\lfloor\frac{m-1}{3}\right\rfloor \leq  2 + 3^{k-4} + 2\left\lfloor\frac{3^{k-2} + 3^{k-3}-1}{3}\right\rfloor
\\
= & 2  + 3^{k-4} + 2(3^{k-2} + 3^{k-3}-1) =  3^{k-2}.
\end{align*}

\item If $y> 3^{k-2} + 3^{k-3}$ then the result holds by part (1). Thus we assume $y \leq 3^{k-2} + 3^{k-3}. $

We have

$$   G(m, y-1) =  \begin{cases}

3^{k-2} , &\text{if $3^{k-2} < m \leq  3^{k-2} + 3^{k-3}$ } \\

3^{k-2} +  \left(\left\lfloor\frac{y-1}{3}\right\rfloor - \left\lfloor\frac{y-2}{3}\right\rfloor \right)(m - 4\cdot 3^{k-3}  )  , &\text{if $3^{k-2} + 3^{k-3} < m \leq  3^{k-2} + 2\cdot 3^{k-3}$ } \\
    
3^{k-2} +   \left(\left\lfloor\frac{y-1}{3}\right\rfloor - \left\lfloor\frac{y-2}{3}\right\rfloor \right)(m - 4\cdot 3^{k-3}  )   , &\text{if $  3^{k-2} + 2\cdot 3^{k-3} \leq m \leq  2\cdot 3^{k-2}$ } \\

m - 3^{k-2}  + \left(\left\lfloor\frac{y-1}{3}\right\rfloor - \left\lfloor\frac{y-2}{3}\right\rfloor \right)(14\cdot 3^{k-3}  -2m)  , &\text{if $ 2\cdot 3^{k-2} \leq m \leq 2\cdot 3^{k-2} + 3^{k-3}$ } \\

m - 3^{k-2} \geq 3^{k-2} + 3^{k-3}  , &\text{if $  2\cdot 3^{k-2} + 3^{k-3} \leq m \leq 2\cdot 3^{k-2} + 2\cdot 3^{k-3}$ } \\

m - 3^{k-2} \geq 3^{k-2} + 2\cdot 3^{k-3}   , &\text{if $  2\cdot 3^{k-2} + 2\cdot 3^{k-3} \leq m \leq 3^{k-1}$ } 

\end{cases}$$

Thus for $3^{k-2} + 3^{k-3} < m \leq  2\cdot 3^{k-2}$  we have

$$   G(m, y-1) =  \begin{cases}

3^{k-2} , &\text{if $y = 1 \mod 3,$ } \\

m - 3^{k-3} \geq 3^{k-2} , &\text{otherwise. } \\

\end{cases}$$

For $ 2\cdot 3^{k-2} \leq m \leq 2\cdot 3^{k-2} + 3^{k-3}$  we have

$$   G(m, y-1) =  \begin{cases}

m - 3^{k-2} , &\text{if $y = 1\mod 3,$ } \\

11\cdot 3^{k-3} - m\geq 3^{k-2} + 3^{k-3} , &\text{otherwise. } \\

\end{cases}$$

\item 

Follows immediately from \cref{nsmol}.

\end{enumerate}
\end{proof}

\begin{lemma}\label{threshold1}
Let $k \geq 4$ and suppose that $A(n',m') = B(n',m')$  for all $n' \leq 3^{k-1}.$ If $1+3^{k-1}< n \leq 3^{k-1} + 3^{k-2},$ and $3^{k-2} \leq m$ then there exists an optimal move $(a,b,c)$ with $3^{k-3}\leq a\leq b\leq c.$ 
\end{lemma}

\begin{proof}
Let $(a,b,c)$ be an optimal move for the game state $(n,m)$ such that $a\leq b\leq c.$ Since $1+3^{k-1}< n$ we have $c > 3^{k-2}.$ 

Suppose $b < 3^{k-3}.$ Thus by \cref{allgaps} $(2) + (3)$ we have

\begin{align*}
S(a,b+1,c-1) - S(a,b,c) = & (A(m,c) - A(m, c-1)) - (A(m,b+1) - A(m, b)) 
\\
= & G(m,c-1) - G(m,b)\geq  3^{k-2} - 3^{k-2} = 0.
\end{align*}

Thus $(a,b+1,c-1)$ is an optimal move. Repeat this procedure until we arrive at an optimal move $(a',b',c')$ with $a'\leq b' \leq c' $ and $3^{k-3} \leq b'.$

Now suppose $a < 3^{k-3}$ and $b\geq 3^{k-3}.$ If $b=c$ then  $(a+1,b-1,c)$ is an optimal move by similar reasoning as above. If $b<c$ then  $(a+1,b,c-1)$ is an optimal move by similar reasoning as above. By repeating these two procedures we arrive at an optimal move $(a'',b'',c'')$ with $ 3^{k-3} \leq a'\leq b' \leq c' .$

\end{proof}

\begin{lemma}\label{threshold2}
Let $k \geq 4$ and suppose that $A(n',m') = B(n',m')$ for all $n' \leq 3^{k-1}.$ If $1+3^{k-1}\leq n \leq 3^{k}$ and $3^{k-2} + 1 < m\leq 3^{k-1}$ then there exists an optimal move $(a,b,c)$ with $3^{k-2}\leq a\leq b\leq c.$ 
\end{lemma}

\begin{proof}

Let $(a,b,c)$ be an optimal move with $a\leq b\leq c$ and  $a < 3^{k-2}.$ By \cref{threshold1} we may assume that $a\geq 3^{k-3}.$ Note that $c>3^{k-2}.$

Suppose $b < 3^{k-2}.$

We first consider the case $1+3^{k-2} \leq m \leq 3^{k-2} + 3^{k-3}.$  By \cref{allgaps} $(1) + (3)$ we have 
$G(m,b) \leq 3^{k-2} \leq  G(m,c-1).$ We can then proceed as in the proof of \cref{threshold1}.

Now consider the case $ 3^{k-2} + 3^{k-3} \leq m \leq 3^{k-1}.$ 

For $3^{k-3} < x < 3^{k-2}, y > 3^{k-2}$  we have

$$  G(m,x) = \begin{cases}

3^{k-2} +  \left(\left\lfloor\frac{x}{3}\right\rfloor - \left\lfloor\frac{x-1}{3}\right\rfloor \right)(m - 4\cdot 3^{k-3}  )  , &\text{if $3^{k-2} + 3^{k-3} < m \leq  3^{k-2} + 2\cdot 3^{k-3}$ } \\
    
3^{k-2} +  \left(\left\lfloor\frac{x}{3}\right\rfloor - \left\lfloor\frac{x-1}{3}\right\rfloor \right)(m - 4\cdot 3^{k-3}  )   , &\text{if $  3^{k-2} + 2\cdot 3^{k-3} \leq m \leq  2\cdot 3^{k-2}$ } \\

m - 3^{k-2}  + \left(\left\lfloor\frac{x}{3}\right\rfloor - \left\lfloor\frac{x-1}{3}\right\rfloor \right)(14\cdot 3^{k-3}  -2m)  , &\text{if $ 2\cdot 3^{k-2} \leq m \leq 2\cdot 3^{k-2} + 3^{k-3}$ } \\

m - 3^{k-2}  + \left(\left\lfloor\frac{x}{2}\right\rfloor - \left\lfloor\frac{x-1}{2}\right\rfloor \right)(14\cdot 3^{k-3}  -2m) \geq m - 3^{k-2}   , &\text{if $  2\cdot 3^{k-2} + 3^{k-3} \leq m \leq 2\cdot 3^{k-2} + 2\cdot 3^{k-3}$ } \\

 3^{k-2} + 2\cdot 3^{k-3}\left(\left\lfloor\frac{x}{2}\right\rfloor - \left\lfloor\frac{x-1}{2}\right\rfloor \right) \leq 3^{k-2} + 2\cdot 3^{k-3}, &\text{if $  2\cdot 3^{k-2} + 2\cdot 3^{k-3} \leq m \leq 3^{k-1}$ } 

\end{cases}$$

$$   G(m, y-1) =  \begin{cases}

3^{k-2} +  \left(\left\lfloor\frac{y-1}{3}\right\rfloor - \left\lfloor\frac{y-2}{3}\right\rfloor \right)(m - 4\cdot 3^{k-3}  )  , &\text{if $3^{k-2} + 3^{k-3} < m \leq  3^{k-2} + 2\cdot 3^{k-3}$ } \\
    
3^{k-2} +   \left(\left\lfloor\frac{y-1}{3}\right\rfloor - \left\lfloor\frac{y-2}{3}\right\rfloor \right)(m - 4\cdot 3^{k-3}  )   , &\text{if $  3^{k-2} + 2\cdot 3^{k-3} \leq m \leq  2\cdot 3^{k-2}$ } \\

m - 3^{k-2}  + \left(\left\lfloor\frac{y-1}{3}\right\rfloor - \left\lfloor\frac{y-2}{3}\right\rfloor \right)(14\cdot 3^{k-3}  -2m)  , &\text{if $ 2\cdot 3^{k-2} \leq m \leq 2\cdot 3^{k-2} + 3^{k-3}$ } \\

m - 3^{k-2}  , &\text{if $  2\cdot 3^{k-2} + 3^{k-3} \leq m \leq 2\cdot 3^{k-2} + 2\cdot 3^{k-3}$ } \\

m - 3^{k-2} \geq 3^{k-2} + 2\cdot 3^{k-3}   , &\text{if $  2\cdot 3^{k-2} + 2\cdot 3^{k-3} \leq m \leq 3^{k-1}$ } 

\end{cases}$$

Thus for $m\geq 2\cdot 3^{k-2} + 3^{k-3} $ we have $G(m,b) \leq G(m,c-1)$ and we proceed as before.

If $ m \leq  2\cdot 3^{k-2} + 3^{k-3} $ and $b \leq 3^{k-2} - 3$  then by the above and \cref{allgaps} $(1)$ we have that $(a+3,b,c-3)$ is an optimal move. Repeating this procedure (and re-ording the tuple to be in order if needed) allows us to assume $b =3^{k-2} - 2$ or  $b =3^{k-2} - 1.$ However in both cases we have $G(m,b) \leq G(m,c-1)$ and we can continue as above.

The case $b \geq 3^{k-2}$ and $a< 3^{k-2}$ follows similarly.
\end{proof}

\begin{lemma}\label{threshold3}
Let $k \geq 4$ and suppose that $A(n',m') = B(n',m')$  for all $n' \leq 3^{k-1}.$ If $3^{k-1} +  3^{k-2} \leq n \leq  3^{k}$ and $3^{k-2} + 1 < m\leq 3^{k-1}$.

Then there exists an optimal move $(a,b,c)$ with $3^{k-2} + 3^{k-3}\leq a\leq b\leq c.$ 
\end{lemma}

\begin{proof}

 By \cref{threshold2} there exists an optimal solution $(a,b,c)$ with $3^{k-2} \leq a \leq b \leq c.$

Suppose $b < 3^{k-2} + 3^{k-3}.$ 

If $b \leq 3^{k-2} + 3^{k-3} -3$ then by the induction tables we have

$$    A(m,b+3) - A(m,b) = \begin{cases}

 3^{k-1}, &\text{if $ 1 + 3^{k-2}  <  m \leq 3^{k-2} + 3^{k-3}  $ } \\

 m + 5 \cdot 3^{k-3}, &\text{if $  3^{k-2}  + 3^{k-3} \leq m \leq 2\cdot 3^{k-2}  + 3^{k-3} $ } \\

 3m - 3^{k-1}  , &\text{if $   2\cdot 3^{k-2} + 3^{k-3} \leq m \leq  3^{k-1} $  } 

\end{cases}$$

Then by \cref{allgaps} we have

$$A(m,c) - A(m,c-3) \geq 3\left( m + 3^{k-3} - 1 -\left\lfloor\frac{m-1}{3}\right\rfloor \right)\geq 2m + 3^{k-2} \geq  A(m,b+3) - A(m,b) .$$

If $b = 3^{k-2} + 3^{k-3} - 2$ or $b = 3^{k-2} + 3^{k-3} - 1$ we have $A(m,c) - A(m,c-1)\geq A(m,b+1) - A(m,b).$

Using these facts we finish the proof as \cref{threshold2}.
\end{proof}

\begin{lemma}\label{optimalform1}
Let $k \geq 4$ and suppose that $A(n',m) = B(n',m)$  for all $n' \leq 3^{k-1}.$ Let $ 3^{k-1} < n \leq 3^{k} $ and $3^{k-2} + 1 < m\leq 3^{k-1}.$ 

Define 

$$\ell = \left\lfloor\frac{n - 3^{k-1}}{3^{k-2}}\right\rfloor.$$

Then there exists an optimal move $(a,b,c)$  for $(m,n)$ with $3^{k-2} + \ell \cdot 3^{k-3} \leq a \leq b \leq c \leq 3^{k-2} + (\ell +1) 3^{k-3}.$ 

Furthermore if $\ell = 0$ then we may assume

$$\left\lfloor\frac{a-1}{3}\right\rfloor+\left\lfloor\frac{b-1}{3}\right\rfloor+\left\lfloor\frac{c-1}{3}\right\rfloor = \left\lfloor\frac{n-1}{3}\right\rfloor - 2.$$

\end{lemma}

\begin{proof}

First suppose we have $1 + 3^{k-1} < n \leq 3^{k-1} + 3^{k-2}. $ Then $\ell = 0.$ By \cref{threshold2} there exists an optimal move $(a,b,c)$ with $3^{k-2} \leq a\leq b \leq c$ and by performing a series of move switches as in \cref{threshold2} we can assume that $c\leq 3^{k-2} + 3^{k-3}.$

If at least two of $(a,b,c)$ are divisible by 3 then 

$$\left\lfloor\frac{a-1}{3}\right\rfloor+\left\lfloor\frac{b-1}{3}\right\rfloor+\left\lfloor\frac{c-1}{3}\right\rfloor = \left\lfloor\frac{n-1}{3}\right\rfloor - 2$$
holds.

Suppose $x$ and $y$ are distinctly indexed elements of $(a,b,c)$ not divisible by $3.$ Note we have $3^{k-2} < x < 3^{k-2} + 3^{k-3}$ and $3^{k-2} < y < 3^{k-2}3^{k-3}$.

If $x \mod 3 = 1$ then by applying the induction formula we have $A(m, y + 1) - A(m, y) \leq A( m, x ) - A(m, x-1). $

On the other hand if $x \mod 3 = 2$ then by applying the induction formula we have $A(m, x + 1) - A(m, x) \leq A( m, y ) - A(m, y-1). $

Thus we can switch to an optimal move with at least one more element divisible by $3$. By repeating this procedure once more if needed we arrive at the desired optimal move.

If instead $ 3^{k-1} + 3^{k-2} \leq n \leq 3^{k-1} $ we first apply \cref{threshold2} to find an optimal move $(a,b,c)$ with $3^{k-2} + 3^{k-3} \leq a\leq b \leq c$ and note that by \cref{allgaps} $(1)$ the moves $(a+1,b-1,c), (a+1,b,c-1),(a,b+1,c-1)$ are all optimal and thus a sequence of move switches shows the desired bounds can be attained.
\end{proof}

\begin{prop}\label{maintheorem1}
Let $k \geq 4$ and suppose that $A(n',m') = B(n',m')$for all $n' \leq 3^{k-1}.$ Then $A( n,m) = B(n,m)$  for all $1+3^{k-1}< n \leq 3^{k}$ and $1+3^{k-2}< m \leq 3^{k-1}.$ 

\end{prop}

\begin{proof}

We break down the proof in to the following cases:

\begin{enumerate}

\item $n = 1 + 3^{k-1}, 1 + 3^{k-2} \leq m \leq 3^{k-1}.$

Let $(a,b,c)$ be the optimal move defined in \cref{optimalform1}.

Then $(a,b,c) = (3^{k-2},3^{k-2},1+3^{k-2}).$

By induction

$$ B(1+3^{k-2}, z) = \begin{cases}

(1 + 2\cdot 3^{k-3})z - 3^{k-3}( 1 +  3^{k-3} ), & \text{if $ 3^{k-3} \leq  z \leq 3^{k-2} + 3^{k-4} $ } \\

  7z\cdot 3^{k-4} + 3^{k-4}( 1 - 13\cdot 3^{k-4} ), & \text{if $ 3^{k-3} + 3^{k-4} \leq z \leq  2\cdot 3^{k-3} $ } \\

(2 + 7\cdot 3^{k-4})z - 3^{k-4}( 11 + 13\cdot 3^{k-4} ), & \text{if $ 2\cdot 3^{k-3} \leq z \leq 2\cdot 3^{k-3} + 3^{k-4} $ } \\

m\cdot 3^{k-2} +  3^{k-3}( 1 - 3^{k-2} ), & \text{if $  2\cdot 3^{k-3} + 3^{k-2}  \leq z \leq 3^{k-2} $ } 

\end{cases}
$$

and so

$$ B(m, 1+3^{k-2}) = \begin{cases}
m\cdot 3^{k-3} + 3^{k-2}( 1 +  3^{k-3} )  , & \text{if $1 + 3^{k-2} \leq m \leq 3^{k-2} + 3^{k-3}$ } \\

(1+2\cdot 3^{k-4})m - 3^{k-3}( 1 - 13\cdot 3^{k-4} )   , & \text{if $ 3^{k-2} + 3^{k-3} \leq m \leq  2\cdot 3^{k-2}$ } \\

(-1 + 2\cdot 3^{k-4})m + 3^{k-3}( 11 + 13\cdot 3^{k-4} )  , & \text{if $  2\cdot 3^{k-2} \leq m \leq 2\cdot 3^{k-2}+ 3^{k-3}$ } \\

m - 3^{k-2}( 1 - 3^{k-2} )  , & \text{if $  2\cdot 3^{k-2} + 3^{k-3} \leq m \leq 3^{k-1}$ } 
\end{cases}
$$

We also have 

$$ A(m, 3^{k-2}) = \begin{cases}

3^{k-3}(m + 3^{k-2}), & \text{if $ 3^{k-2}  \leq  m \leq 3^{k-2} + 3^{k-3} $ } \\

2m\cdot 3^{k-4} + 13\cdot 3^{2k-7} , & \text{if $ 3^{k-2} + 3^{k-3} \leq m \leq  2\cdot 3^{k-2} $ } \\

 2m\cdot 3^{k-4} + 13\cdot 3^{2k-7} , & \text{if $ 2\cdot 3^{k-2} \leq m \leq 2\cdot 3^{k-2} + 3^{k-3} $ } \\

 m\cdot 3^{k-1} , & \text{if $  2\cdot 3^{k-2} +  3^{k-3}  \leq   m \leq 3^{k-1} $ } 

\end{cases}
$$

Therefore

$$ B(1+3^{k-1}, m) = \begin{cases}

(1 + 2\cdot 3^{k-2})m - 3^{k-2}( 1 +  3^{k-2} ), & \text{if $ 3^{k-2} \leq  m \leq 3^{k-2} + 3^{k-3} $ } \\

  7m\cdot 3^{k-3} + 3^{k-3}( 1 - 13\cdot 3^{k-3} ), & \text{if $ 3^{k-2} + 3^{k-3} \leq m \leq  2\cdot 3^{k-2} $ } \\

(2 + 7\cdot 3^{k-3})m - 3^{k-3}( 11 + 13\cdot 3^{k-3} ), & \text{if $ 2\cdot 3^{k-2} \leq m \leq 2\cdot 3^{k-2} + 3^{k-3} $ } \\

m\cdot 3^{k-1} +  3^{k-2}( 1 - 3^{k-1} ), & \text{if $  2\cdot 3^{k-2} + 3^{k-3}  \leq m \leq  3^{k-1} $ } 

\end{cases}
$$

as desired. 

\item $1 + 3^{k-1} \leq n\leq 3^{k},$  $m = 1 + 3^{k-2}.$

$$ B(n, 1+3^{k-2}) = \begin{cases}
n\cdot 3^{k-2} + 3^{k-2}( 1 +  3^{k-3} )  , & \text{if $1 + 3^{k-1} \leq n \leq 3^{k-1} + 3^{k-2}$ } \\

(1+2\cdot 3^{k-4})n - 3^{k-3}( 1 - 13\cdot 3^{k-4} )   , & \text{if $ 3^{k-1} + 3^{k-2} \leq n \leq  2\cdot 3^{k-1}$ } \\

(-1 + 2\cdot 3^{k-4})n + 3^{k-3}( 11 + 13\cdot 3^{k-4} )  , & \text{if $  2\cdot 3^{k-1} \leq n \leq 2\cdot 3^{k-1}+ 3^{k-2}$ } \\

n - 3^{k-4}( 1 - 3^{k-4} )  , & \text{if $  2\cdot 3^{k-1} + 3^{k-2} \leq n \leq 3^{k}$ } 
\end{cases}
$$

\item $ 1 + 3^{k-1} < n \leq  3^{k-1} + 3^{k-2},$  $ 1 + 3^{k-2} < m \leq 3^{k-2} + 3^{k-3}.$

Let $(a,b,c)$ be the optimal move defined in \cref{optimalform1}.

We have $3^{k-2} \leq a \leq b \leq c \leq 3^{k-2}  + 3^{k-3}.$ Therefore 

\begin{align*}
A(n,m) & = S(a,b,c,m) 
\\
& = nm - A(m,a) - A(m,b) - A(m,c)
\\
& = nm - n\cdot 3^{k-2} - 3^{k-2}(m - 2\cdot 3^{k-2} )
\\
& = (n - 3^{k-2})m - 3^{k-2}(n - 2\cdot 3^{k-2})
\\
& = B(n,m).
\end{align*}

\item $ 1 + 3^{k-1} < n \leq  3^{k-1} + 3^{k-2},$  $ 3^{k-2} + 3^{k-3} \leq m \leq  2\cdot 3^{k-2}.$

Let $(a,b,c)$ be the optimal move defined in \cref{optimalform1}.

We have $3^{k-2} \leq a \leq b \leq c \leq 3^{k-2}  + 3^{k-3}$ and

$$\left\lfloor\frac{a-1}{3}\right\rfloor+\left\lfloor\frac{b-1}{3}\right\rfloor+\left\lfloor\frac{c-1}{3}\right\rfloor = \left\lfloor\frac{n-1}{3}\right\rfloor - 2.$$

Therefore 

\begin{align*}
A(n,m) = & S(a,b,c,m) 
\\
= &  nm - A(m,a) - A(m,b) - A(m,c)
\\
 = & nm - n\cdot 3^{k-2} - \left(\left\lfloor\frac{n-1}{3}\right\rfloor - 2\right)(m -4\cdot 3^{k-3} ) 
\\
& - 3m + 3^{k-3}( m  + 12 + 2\cdot 3^{k-3} )
\\
 = & \left(n-1 + 3^{k-3} - \left\lfloor\frac{n-1}{3}\right\rfloor \right)m + 4\cdot 3^{k-3}\left\lfloor\frac{n-1}{3}\right\rfloor 
\\ 
& + 3^{k-3}( -3n + 4  + 2\cdot 3^{k-3})
 \\
= &  B(n,m).
\end{align*}

\item $ 1 + 3^{k-1} < n \leq  3^{k-1} + 3^{k-2},$  $  2\cdot 3^{k-2} \leq m \leq  2\cdot 3^{k-2} + 3^{k-3} .$

Let $(a,b,c)$ be the optimal move defined in \cref{optimalform1}.

We have $3^{k-2} \leq a \leq b \leq c \leq 3^{k-2}  + 3^{k-3}$ and 

$$\left\lfloor\frac{a-1}{3}\right\rfloor+\left\lfloor\frac{b-1}{3}\right\rfloor+\left\lfloor\frac{c-1}{3}\right\rfloor = \left\lfloor\frac{n-1}{3}\right\rfloor - 2.$$

Therefore 

\begin{align*}
A(n,m) = & S(a,b,c,m) 
\\
= &  nm - A(m,a) - A(m,b) - A(m,c)
\\
 = & nm - nm + n\cdot 3^{k-2} - \left(\left\lfloor\frac{n-1}{3}\right\rfloor - 2\right)(14\cdot 3^{k-3} - 2m )
\\
& + 6m - 3^{k-3} \left( - m + 42 - 2\cdot 3^{k-3}  \right)
\\
 = & \left(2 + 3^{k-3} + 2\left\lfloor\frac{n-1}{3}\right\rfloor \right)m  - 14\cdot 3^{k-3}\left\lfloor\frac{n-1}{3}\right\rfloor
\\ 
& + 3^{k-3}( 3n  - 14  + 2\cdot 3^{k-3})
 \\
= &  B(n,m).
\end{align*}

\item $ 1 + 3^{k-1} < n \leq  3^{k-1} + 3^{k-2},$  $   2\cdot 3^{k-2} + 3^{k-3} \leq m \leq 3^{k-1}.$

Let $(a,b,c)$ be the optimal move defined in \cref{optimalform1}.

We have $3^{k-2} \leq a \leq b \leq c \leq 3^{k-2}  + 3^{k-3}.$ Therefore 

\begin{align*}
A(n,m) & = S(a,b,c,m) 
\\
& = nm - A(m,a) - A(m,b) - A(m,c)
\\
& = nm - n(m - 3^{k-2}  ) -  3^{k-1}( 2\cdot 3^{k-2} - m  )
\\
& = m\cdot 3^{k-1} + 3^{k-2}( n - 2\cdot 3^{k-1})
\\
& = B(n,m).
\end{align*}

\item $ 3^{k-1} + 3^{k-2} < n \leq  3^{k-1} + 2\cdot 3^{k-2},$  $1 + 3^{k-2} <  m \leq 3^{k-1}.$

Let $(a,b,c)$ be the optimal move defined in \cref{optimalform1}.

We have $3^{k-2} + 3^{k-3} \leq a \leq b \leq c \leq 3^{k-2}  + 2\cdot 3^{k-3}.$ Therefore 

\begin{align*}
A(n,m) = &  S(a,b,c,m) 
\\
= &  nm - A(m,a) - A(m,b) - A(m,c)
\\
= &  nm - n\left(m + 3^{k-3} - 1 -\left\lfloor\frac{m-1}{3}\right\rfloor \right)
\\
& - 3^{k-2}\left( - 3m + 4 + 2\cdot 3^{k-3} + 4\left\lfloor\frac{m-1}{3}\right\rfloor \right)
\\
= &   m\cdot 3^{k-1} +  \left\lfloor\frac{m-1}{3}\right\rfloor(n - 4\cdot 3^{k-2}  ) +  n  - 3^{k-3}\left( n + 12 +  2\cdot 3^{k-2} \right)
\\
= &  B(n,m).
\end{align*}

\item $   3^{k-1} + 2\cdot 3^{k-2} \leq n \leq  2\cdot 3^{k-1},$  $1+ 3^{k-2} <  m \leq 3^{k-1}.$

Let $(a,b,c)$ be the optimal move defined in \cref{optimalform1}.

We have $3^{k-2} + 2\cdot 3^{k-3} \leq a \leq b \leq c \leq  2\cdot 3^{k-2}.$ Therefore 

\begin{align*}
A(n,m) = &  S(a,b,c,m) 
\\
= &  nm - A(m,a) - A(m,b) - A(m,c)
\\
= &  nm - n\left(m + 3^{k-3} - 1 - \left\lfloor\frac{m-1}{3}\right\rfloor \right)
\\
& - 3^{k-2}\left( - 3m + 4 + 2\cdot 3^{k-3} + 4\left\lfloor\frac{m-1}{3}\right\rfloor \right)
\\
= &   m\cdot 3^{k-1} +  \left\lfloor\frac{m-1}{3}\right\rfloor(n - 4\cdot 3^{k-2}  ) +  n  - 3^{k-3}\left( n + 12 + 2\cdot 3^{k-2} \right) 
\\
= &  B(n,m).
\end{align*}

\item $  2\cdot 3^{k-1} \leq n \leq 2\cdot 3^{k-1} + 3^{k-2},$  $1+ 3^{k-2} <  m \leq 3^{k-1}.$

Let $(a,b,c)$ be the optimal move defined in \cref{optimalform1}.

We have $ 2\cdot 3^{k-2}  \leq a \leq b \leq c \leq 3^{k-1}.$ Therefore

\begin{align*}
A(n,m) = &  S(a,b,c,m) 
\\
= &  nm - A(m,a) - A(m,b) - A(m,c)
\\
= &  nm - n\left(2 + 3^{k-3} + 2\left\lfloor\frac{m-1}{3}\right\rfloor  \right)
\\
& - 3^{k-2}\left( 3m - 14 + 2\cdot 3^{k-3} - 14\left\lfloor\frac{m-1}{3}\right\rfloor    \right) 
\\
= &   (n - 3^{k-1} )m +\left\lfloor\frac{m-1}{3}\right\rfloor(14\cdot 3^{k-2}  -2n) - 2n - 3^{k-3} \left( n - 42 + 2\cdot 3^{k-2}  \right)
\\
= &  B(n,m).
\end{align*}

\item $ 2\cdot 3^{k-1} + 3^{k-2} \leq n \leq 2\cdot 3^{k-1} + 2\cdot 3^{k-2} ,$  $1+ 3^{k-2} <  m \leq 3^{k-1}.$

\begin{align*}
A(n,m) = &  S(a,b,c,m) 
\\
= &  nm - A(m,a) - A(m,b) - A(m,c)
\\
= &  nm - n\left(1 + 2\left\lfloor\frac{m-1}{2}\right\rfloor\right)
\\
& - 3^{k-2} \left( 3m - 7 +  3^{k-1} - 14\left\lfloor\frac{m-1}{2}\right\rfloor \right)
\\
= &   (n - 3^{k-1} )m  + \left\lfloor\frac{m-1}{2}\right\rfloor( 14\cdot 3^{k-2} - 2n  ) - n + 3^{k-2} \left( 7 - 3^{k-1}\right)
\\
= &  B(n,m).
\end{align*}

\item $2\cdot 3^{k-1} + 2\cdot 3^{k-2} \leq n \leq 3^{k},$  $1 + 3^{k-2} <  m \leq 3^{k-1}.$

Let $(a,b,c)$ be the optimal move defined in \cref{optimalform1}.

We have $2\cdot 3^{k-2} + 2\cdot 3^{k-3}  \leq a \leq b \leq c \leq 3^{k-1}.$ Therefore 

\begin{align*}
A(n,m) & = S(a,b,c,m) 
\\
& = nm - A(m,a) - A(m,b) - A(m,c)
\\
& = nm - nm - 2\cdot 3^{k-2}\left\lfloor\frac{m-1}{2}\right\rfloor + 3^{k-2}\left( 5m - 1 - 3^{k-2} \right)
\\
& =   m\cdot 3^{k-1} + 2\cdot 3^{k-2}\left\lfloor\frac{m}{2}\right\rfloor + 3^{k-2} - 3^{2k-3} 
\\
& = B(n,m).
\end{align*}

\end{enumerate}
\end{proof}

\begin{lemma}\label{inductionlemma}
Let $k \geq 4$ and suppose that $A(n',m') = B(n',m')$for all $n' \leq 3^{k-1}.$ Then

\begin{enumerate}

\item $ 1 + 3^{k-1} < m \leq  3^{k-1} + 3^{k-2},$  $ 3^{k-2} < x \leq 3^{k-2} + 3^{k-3}.$

$$ A(m, x) = (m - 3^{k-2})x - 3^{k-2}(m - 2\cdot 3^{k-2}). $$

\item $ 1 + 3^{k-1} < m \leq  3^{k-1} + 3^{k-2},$  $ 3^{k-2} + 3^{k-3} \leq x  \leq 2\cdot 3^{k-2}.$

$$A(m,x) = \left(m-1 + 3^{k-3} - \left\lfloor\frac{m-1}{3}\right\rfloor \right)x + 4\cdot 3^{k-3}\left\lfloor\frac{m-1}{3}\right\rfloor + 3^{k-3}( -3m + 4  + 2\cdot 3^{k-3}).$$

\item $ 1 + 3^{k-1} < m \leq  3^{k-1} + 3^{k-2},$  $  2\cdot 3^{k-2} \leq x \leq  2\cdot 3^{k-2} + 3^{k-3} .$

$$A(m,x) = \left(2 + 3^{k-3} + 2\left\lfloor\frac{m-1}{3}\right\rfloor \right)x  - 14\cdot 3^{k-3}\left\lfloor\frac{m-1}{3}\right\rfloor + 3^{k-3}( 3m  - 14  + 2\cdot 3^{k-3}).$$

\item $ 1 + 3^{k-1} < m \leq  3^{k-1} + 3^{k-2},$  $   2\cdot 3^{k-2} + 3^{k-3} \leq x  \leq 3^{k-1}.$

$$A(m,x) = x\cdot 3^{k-1} + 3^{k-2}( m - 2\cdot 3^{k-1}). $$

\end{enumerate}

\end{lemma}

\begin{proof}
Follows directly from \cref{maintheorem1}.
\end{proof}

\begin{lemma}\label{optimalform2}
Let $k \geq 4$ and suppose that $A(n',m) = B(n',m)$  for all $n' \leq 3^{k-1}.$ Let $ 3^{k-1} < n \leq 3^{k} $ and $3^{k-1} + 1 < m\leq 3^{k-1} + 3^{k-2}.$ 

Define 

$$\ell = \left\lfloor\frac{n - 3^{k-1}}{3^{k-2}}\right\rfloor.$$

Then there exists an optimal move $(a,b,c)$  for $(m,n)$ with $3^{k-2} + \ell \cdot 3^{k-3} \leq a \leq b \leq c \leq 3^{k-2} + (\ell +1) 3^{k-3}.$ 

Furthermore if $\ell = 0$ then we may assume

$$\left\lfloor\frac{a-1}{3}\right\rfloor+\left\lfloor\frac{b-1}{3}\right\rfloor+\left\lfloor\frac{c-1}{3}\right\rfloor = \left\lfloor\frac{n-1}{3}\right\rfloor - 2.$$

\end{lemma}

\begin{proof}

By \cref{inductionlemma} we have

$$
A(m,x) - A(m,x-1) = 
\begin{cases}
\Phi(x)-\Phi(x-1) \leq  5\cdot 3^{k-3}, &\text{if $ x \leq 3^{k-2}$ } \\

m - 3^{k-2}, &\text{if  $ 3^{k-2} < x \leq 3^{k-2} + 3^{k-3}$ } \\

m-1 + 3^{k-3} - \left\lfloor\frac{m-1}{3}\right\rfloor , &\text{if  $ 3^{k-2} + 3^{k-3} \leq x  \leq 2\cdot 3^{k-2}$ }\\

2 + 3^{k-3} + 2\left\lfloor\frac{m-1}{3}\right\rfloor, &\text{if $  2\cdot 3^{k-2} \leq x \leq  2\cdot 3^{k-2} + 3^{k-3}$ } \\

3^{k-1}, &\text{if  $  2\cdot 3^{k-2} + 3^{k-3} \leq x  \leq 3^{k-1}$}

\end{cases}
$$

We have

$$3^{k-2}\leq m - 3^{k-2} \leq m-1 + 3^{k-3} - \left\lfloor\frac{m-1}{3}\right\rfloor \leq 2 + 3^{k-3} + 2\left\lfloor\frac{m-1}{3}\right\rfloor \leq 3^{k-1}.$$

Therefore for $x\leq y \leq 3^{k-1}$ we have $G(m,x) \leq G(m,y)$ and we can proceed as in \cref{optimalform1}.
\end{proof}

\begin{prop}\label{maintheorem2}
Let $k \geq 4$ and suppose that $A(n',m') = B(n',m')$for all $n' \leq 3^{k-1}.$ Then $A( n,m) = B(n,m)$  for all $1+3^{k-1}\leq n \leq 3^{k}$ and $1+3^{k-1}< m < 3^{k-1} + 3^{k-2}.$ 

\end{prop}

\begin{proof}
We break down the proof in to the following cases:

\begin{enumerate}

\item $ 1 + 3^{k-1} < n \leq 3^{k-1} + 3^{k-2},$  $1+ 3^{k-1} \leq  m < 3^{k-1} + 3^{k-2}.$

Let $(a,b,c)$ be the optimal move defined in \cref{optimalform2}.
We have $1+ 3^{k-2} < a \leq b \leq c \leq 3^{k-2}  + 3^{k-3}.$ Therefore 

\begin{align*}
A(n,m) = &  S(a,b,c,m) 
\\
= &  nm - A(m,a) - A(m,b) - A(m,c)
\\
= &  nm - n(m - 3^{k-2}) +  3^{k-1}(m - 2\cdot 3^{k-2})
\\
= &   m\cdot 3^{k-1} + 3^{k-2}(n - 2\cdot 3^{k-1} )
\\
= &  B(n,m).
\end{align*}

\item $ 3^{k-1} + 3^{k-2} \leq n \leq 3^{k-1} + 2 \cdot 3^{k-2},$  $1+ 3^{k-1} \leq  m < 3^{k-1} + 3^{k-2}.$

Let $(a,b,c)$ be the optimal move defined in \cref{optimalform2}.
We have $ 3^{k-2} + 3^{k-3} \leq a \leq b \leq c \leq 3^{k-2}  + 2\cdot 3^{k-3}.$ Therefore 

\begin{align*}
A(n,m) = &  S(a,b,c,m) 
\\
= &  nm - A(m,a) - A(m,b) - A(m,c)
\\
= &  nm - \left(m-1 + 3^{k-3} - \left\lfloor\frac{m-1}{3}\right\rfloor \right)n -  4\cdot 3^{k-2}\left\lfloor\frac{m-1}{3}\right\rfloor + 3^{k-3}( -3m + 4  + 2\cdot 3^{k-3})
\\
= &  m\cdot 3^{k-1}+ \left\lfloor\frac{m-1}{3}\right\rfloor( n -4\cdot 3^{k-2}   )  +  n  - 3^{k-3}\left( n  + 12 + 2\cdot 3^{k-2} \right)
\\
= &  B(n,m).
\end{align*}

\item $ 3^{k-1} + 2 \cdot 3^{k-2} \leq n \leq 2 \cdot 3^{k-1},$  $1+3^{k-1} \leq  m < 3^{k-1} + 3^{k-2}.$

Let $(a,b,c)$ be the optimal move defined in \cref{optimalform2}.
We have $ 3^{k-2} + 2\cdot 3^{k-3} \leq a \leq b \leq c \leq 2\cdot 3^{k-2}.$ Therefore

\begin{align*}
A(n,m) = &  S(a,b,c,m) 
\\
= &  nm - A(m,a) - A(m,b) - A(m,c)
\\
= &  nm - \left(m-1 + 3^{k-3} - \left\lfloor\frac{m-1}{3}\right\rfloor \right)n -  4\cdot 3^{k-2}\left\lfloor\frac{m-1}{3}\right\rfloor + 3^{k-3}( -3m + 4  + 2\cdot 3^{k-3})
\\
= &  m\cdot 3^{k-1}+ \left\lfloor\frac{m-1}{3}\right\rfloor( n -4\cdot 3^{k-2}   )  +  n  - 3^{k-3}\left( n  + 12 + 2\cdot 3^{k-2} \right)
\\
= &  B(n,m).
\end{align*}

\item $  2 \cdot 3^{k-1} \leq n \leq  2 \cdot 3^{k-1} + 3^{k-2},$  $1+3^{k-1} \leq  m < 3^{k-1} + 3^{k-2}.$

Let $(a,b,c)$ be the optimal move defined in \cref{optimalform2}.
We have $  2\cdot 3^{k-2} \leq a \leq b \leq c \leq 2\cdot 3^{k-2} + 3^{k-3}.$ Therefore 

\begin{align*}
A(n,m) = &  S(a,b,c,m) 
\\
= &  nm - A(m,a) - A(m,b) - A(m,c)
\\
= &  nm - \left(2 + 3^{k-3} + 2\left\lfloor\frac{m-1}{3}\right\rfloor \right)n  + 14\cdot 3^{k-2}\left\lfloor\frac{m-1}{3}\right\rfloor + 3^{k-3}( 3m  - 14  + 2\cdot 3^{k-3})
\\
= & (n -  3^{k-1}  )m  + \left\lfloor\frac{m-1}{3}\right\rfloor (14\cdot 3^{k-2} - 2n )  - 2n + 3^{k-3} \left( - n + 42 - 2\cdot 3^{k-2}  \right)
\\
= &  B(n,m).
\end{align*}

\item $  2 \cdot 3^{k-1} + 3^{k-2} \leq n \leq 2 \cdot 3^{k-1} + 2\cdot 3^{k-2} ,$  $1+3^{k-1} \leq  m < 3^{k-1} + 3^{k-2}.$

Let $(a,b,c)$ be the optimal move defined in \cref{optimalform2}.
We have $  2\cdot 3^{k-2} +3^{k-3} \leq a \leq b \leq c \leq 2\cdot 3^{k-2} + 2\cdot 3^{k-3}.$ Therefore 

\begin{align*}
A(n,m) = &  S(a,b,c,m) 
\\
= &  nm - A(m,a) - A(m,b) - A(m,c)
\\
= &  nm - 3^{k-1}\cdot n - 3^{k-1}( m - 2\cdot 3^{k-1})
\\
= & (n - 3^{k-1}  )m + 3^{k-1}( 2\cdot 3^{k-1} - n  )
\\
= &  B(n,m).
\end{align*}

\item $ 2 \cdot 3^{k-1} + 2\cdot 3^{k-2} \leq n \leq 3^{k},$  $1+3^{k-1} \leq  m < 3^{k-1} + 3^{k-2}.$

Let $(a,b,c)$ be the optimal move defined in \cref{optimalform2}.
We have $  2\cdot 3^{k-2} + 2\cdot 3^{k-3} \leq a \leq b \leq c \leq 3^{k-1}.$ Therefore 

\begin{align*}
A(n,m) = &  S(a,b,c,m) 
\\
= &  nm - A(m,a) - A(m,b) - A(m,c)
\\
= &  nm - 3^{k-1}\cdot n - 3^{k-1}( m - 2\cdot 3^{k-1})
\\
= & (n - 3^{k-1}  )m + 3^{k-1}( 2\cdot 3^{k-1} - n  )
\\
= &  B(n,m).
\end{align*}
\end{enumerate}
\end{proof}

\begin{lemma}\label{inductionlemma2}
Let $k \geq 4$ and suppose that $A(n',m') = B(n',m')$for all $n' \leq 3^{k-1}.$ Then

\begin{enumerate}

\item $ 3^{k-1} + 3^{k-2} \leq m \leq   2\cdot 3^{k-1},$  $3^{k-2} \leq  x \leq 3^{k-1}.$

$$A(m, 3^{k-2}) = 3^{2k - 4}.$$

$$A(m, 1+3^{k-2}) = 3^{2k - 4} + m - 3^{k-2}.$$

$$A(m, 2+3^{k-2}) = 3^{2k - 4} + m + 2\cdot 3^{k-2}.$$

If $ x \geq 1 + 3^{k-2}$ then
$$A(m,x) = x\cdot 3^{k-1} +  \left\lfloor\frac{x-1}{3}\right\rfloor(m - 4\cdot 3^{k-2}  ) + m - 3^{k-3} \left(  m + 12  + 2\cdot 3^{k-2} \right) $$

$$A(m,3^{k-1} - 1) = 13\cdot 3^{2k-5} - 3^{k-1}  + 2m\cdot 3^{k-3}. $$

$$A(m,3^{k-1}) = 2m\cdot 3^{k-3} + 13\cdot 3^{2k-5}. $$

\item $  2\cdot 3^{k-1} \leq m \leq 2\cdot 3^{k-1} + 3^{k-2},$  $3^{k-2} \leq  x \leq 3^{k-1}.$

$$A(m, 3^{k-2}) = 3^{2k - 4}.$$

$$A(m, 1+3^{k-2}) =  5\cdot 3^{k-2}+  3^{2k - 4}. $$

$$A(m, 2+3^{k-2}) =  8\cdot 3^{k-2}+  3^{2k - 4}. $$

If $ x > 1 + 3^{k-2}$ then
$$A(m,x) = (m - 3^{k-1} )x +\left\lfloor\frac{x-1}{3}\right\rfloor(14\cdot 3^{k-2}  -2m) - 2m - 3^{k-3} \left( m - 42  + 2\cdot 3^{k-2} \right). $$

$$A(m,3^{k-1} - 1) = 2m\cdot 3^{k-3} - m + 3^{k-1} + 13\cdot 3^{2k-5}     $$

$$A(m,3^{k-1}) = 2m\cdot 3^{k-3} + 13\cdot 3^{2k-5}     $$

\item $ 2\cdot 3^{k-1} + 3^{k-2} \leq m \leq 2\cdot 3^{k-1} + 2\cdot 3^{k-2} ,$  $ 3^{k-2} \leq  x \leq 3^{k-1}.$

$$A(m, 3^{k-2}) = 3^{2k - 4}.$$

$$A(m, 1+3^{k-2}) =  5\cdot 3^{k-2}+  3^{2k - 4}. $$

$$A(m, 2+3^{k-2}) =  8\cdot 3^{k-2}+  3^{2k - 4}. $$

If $ x > 1 + 3^{k-2}$ then
$$A(m,x) = (m - 3^{k-1} )x  + \left\lfloor\frac{x-1}{2}\right\rfloor( 14\cdot 3^{k-2} - 2m  ) - m + 3^{k-2} \left( 7 - 3^{k-1}\right).  $$

$$A(m, 3^{k-1} - 1) = m  - 11\cdot 3^{k-2} + 3^{2k-2}  .  $$

$$A(m, 3^{k-1}) =  3^{2k-2}  .  $$

\item $2\cdot 3^{k-1} + 2\cdot 3^{k-2} \leq m \leq 3^{k},$  $3^{k-2} \leq  x \leq 3^{k-1}.$

$$A(m, 3^{k-2}) = 3^{2k - 4}.$$

$$A(m, 1+3^{k-2}) =  5\cdot 3^{k-2}+  3^{2k - 4}. $$

$$A(m, 2+3^{k-2}) =  8\cdot 3^{k-2}+  3^{2k - 4}. $$

If $ x \geq 1 + 3^{k-2}$ then

$$A(m,x) = x\cdot 3^{k-1} + 2\cdot 3^{k-2}\left\lfloor\frac{x}{2}\right\rfloor + 3^{k-2} - 3^{2k-3} . $$

$$A(m, 3^{k-1} - 1) =  3^{2k-2} - 4\cdot 3^{k-2} . $$

$$A(m, 3^{k-1}) =  3^{2k-2} - 3^{k-2} . $$

\end{enumerate}

\end{lemma}

\begin{proof}
Follows directly from \cref{maintheorem1}.
\end{proof}

\begin{lemma}\label{optimalform3}
Let $k \geq 4$ and suppose that $A(n',m') = B(n',m')$for all $n' \leq 3^{k-1}.$ Let $1 + 3^{k-1} \leq n \leq 3^{k-1} + 2\cdot 3^{k-2}$ and $3^{k-1} + 3^{k-2} \leq m. $ Then $(3^{k-2}, 3^{k-2}, n - 2\cdot 3^{k-2})$ is an optimal move for $(n,m).$

\end{lemma}

\begin{proof}

Suppose $(a,b,c)$ is an optimal move with $b < 3^{k-2}.$ Then by \cref{philemma} and \cref{inductionlemma2} we have 

$$A(m, b+1) - A(m, b) \leq 5 \cdot 3^{k-3} \leq 3^{k-1} \leq A(m,c) - A(m,c-1).$$

Thus we can proceed as in \cref{threshold1} and therefore assume we have $3^{k-2} \leq a \leq b\leq c.$

We now break down the proof in to the following cases

\begin{enumerate}

\item $ 3^{k-1} + 3^{k-2} \leq m \leq   2\cdot 3^{k-1}.$ 

Suppose $a\geq 3^{k-2} + 3.$ Then $c\leq 3^{k-1}-3$ and $A(m,c+3) - A(m,c) =   A(m,a)-A(m,a-3).$ Therefore $(a-3,b,c+3)$ is an optimal move.

Suppose $a = 3^{k-2} + 2.$ Then
$$A(m,a)-A(m,a-2) = m + 2\cdot 3^{k-2} \geq 2\cdot 3^{k-1} $$
and

$$
A(m,c+2) - A(m,c) = 
\begin{cases}
2\cdot 3^{k-1}, &\text{if $ c = 1 \mod 3$ } \\

m + 2\cdot 3^{k-2}, &\text{otherwise } 

\end{cases}
$$

Therefore $A(m,c+2) - A(m,c) \leq   A(m,a)-A(m,a-2)$ and $(a-2,b,c+2)$ is an optimal move.

Suppose $a = 3^{k-2} + 1.$ Then $A(m,a)-A(m,a-1) = m - 3^{k-2 } \geq A(c+1,m) - A(c,m) $ and thus $(a-1,b,c+1)$ is an optimal move.

Thus we may assume that $a = 3^{k-2}.$ We may also assume that $b=3^{k-2}$ by similar reasoning.

\item $  2\cdot 3^{k-1} \leq m \leq 2\cdot 3^{k-1} + 3^{k-2},$  $3^{k-2} \leq  x \leq 3^{k-1}.$

Suppose $a\geq 3^{k-2} + 5.$ Then $c\leq 3^{k-1}-5$ and $A(m,c+3) - A(m,c) =   A(m,a)-A(m,a-3).$ 

Suppose $a = 3^{k-2} + 4.$ Then 
$$A(c+3,m) - A(c, m) = m + 5\cdot 3^{k-2}\leq  3^{k} = A(a,m) - A(a-3, m). $$

Suppose $a = 3^{k-2} + 3.$ Then 
$$A(c+3,m) - A(c, m) = m + 5\cdot 3^{k-2}\leq 11\cdot 3^{k-2} = A(a,m) - A(a-3, m). $$

Suppose $a = 3^{k-2} + 2.$ Then 
$$
A(m,c+2) - A(m,c) = 
\begin{cases}
2m - 2\cdot 3^{k-1}\leq 8\cdot 3^{k-2}, &\text{if $ c = 1 \mod 3$ } \\

8\cdot 3^{k-2}, &\text{otherwise } 

\end{cases}
$$
Thus 
$$A(c+2,m) - A(c, m) \leq 8\cdot 3^{k-2} = A(a,m) - A(a-2, m). $$

Suppose $a = 3^{k-2} + 1.$ Then 
$$
A(m,c+1) - A(m,c) = 
\begin{cases}
11\cdot 3^{k-2} - m \leq 5\cdot 3^{k-2}, &\text{if $ c = 0 \mod 3$ } \\

m -  3^{k-1}\leq 4\cdot 3^{k-2}, &\text{otherwise } 

\end{cases}
$$

Thus we may assume that $a = 3^{k-2}.$ We may also assume that $b=3^{k-2}$ by similar reasoning.

\item $ 2\cdot 3^{k-1} + 3^{k-2} \leq m \leq 2\cdot 3^{k-1} + 2\cdot 3^{k-2} ,$  $ 3^{k-2} \leq  x \leq 3^{k-1}.$

Suppose $a\geq 3^{k-2} + 4.$ Then $c\leq 3^{k-1}-4$ and $A(m,c+2) - A(m,c) =   A(m,a)-A(m,a-2).$ 

Suppose $a = 3^{k-2} + 3.$ Then
$$
A(m,c+3) - A(m,c) = 
\begin{cases}
m + 5\cdot 3^{k-2}, &\text{if $ c = 1 \mod 2$ } \\

19\cdot 3^{k-2} - m\leq m + 5\cdot 3^{k-2}, &\text{otherwise } 

\end{cases}
$$
Thus 
$$A(m,c+3) - A(m,c) \leq m + 5\cdot 3^{k-2} = A(a,m) - A(a-3, m). $$

Suppose $a = 3^{k-2} + 2.$ Then
$$A(m,c+2) - A(m,c) \leq 8\cdot 3^{k-2} = A(a,m) - A(a-2, m). $$

Suppose $a = 3^{k-2} + 1.$ Then

$$
A(m,c+1) - A(m,c) = 
\begin{cases}
m - 3^{k-1}\leq 5\cdot 3^{k-2}, &\text{if $ c = 1 \mod 2$ } \\

11\cdot 3^{k-2} - m\leq 4\cdot 3^{k-2}, &\text{otherwise } 

\end{cases}
$$

Thus
$$A(m,c+1) - A(m,c) \leq 5\cdot 3^{k-2} = A(a,m) - A(a-1, m). $$

Thus we may assume that $a = 3^{k-2}.$ We may also assume that $b=3^{k-2}$ by similar reasoning.

\item $2\cdot 3^{k-1} + 2\cdot 3^{k-2} \leq m \leq 3^{k},$  $3^{k-2} \leq  x \leq 3^{k-1}.$

Suppose $a\geq 3^{k-2} + 2.$ Then $c\leq 3^{k-1}-2$ and so
$$A(m,c+2) - A(m,c) =   A(m,a)-A(m,a-2).$$ 

Suppose $a\geq 3^{k-2} + 1.$ Then 
$$A(m,c+1) - A(m,c) \leq 5\cdot 3^{k-2} =  A(m,a)-A(m,a-1).$$

Thus we may assume that $a = 3^{k-2}.$ We may also assume that $b=3^{k-2}$ by similar reasoning.

\end{enumerate}
    
\end{proof}

\begin{prop}\label{maintheorem3}
Let $k \geq 4$ and suppose that $A(n',m') = B(n',m')$for all $n' \leq 3^{k-1}.$ Then $A( n,m) = B(n,m)$  for all $1+3^{k-1}< n \leq 3^{k-1} + 2\cdot 3^{k-2} $ and $ 3^{k-1} + 3^{k-2} \leq m \leq 3^{k}.$ 

\end{prop}

\begin{proof}

By \cref{optimalform3} the move $(3^{k-2}, 3^{k-2}, n - 2\cdot 3^{k-2})$ is optimal for $(n,m).$

We now break down the proof in to the following cases:

\begin{enumerate}

\item  $3^{k-1} +  3^{k-2} \leq m \leq 2\cdot 3^{k-1}  .$

\begin{align*}
A(n,m) = &  S(3^{k-2}, 3^{k-2}, n - 2\cdot 3^{k-2},m) 
\\
= &  nm - 2A(m,3^{k-2})- A(m,n - 2\cdot 3^{k-2})
\\
= &  nm - 2\cdot 3^{2k-4} - 3^{k-1}(n - 2\cdot 3^{k-2}) - m +  3^{k-3} \left(  m + 12  + 2\cdot 3^{k-2} \right)
\\
& - \left\lfloor\frac{n -1}{3}\right\rfloor(m - 4\cdot 3^{k-2}  )  +  2\cdot 3^{k-3}(m - 4\cdot 3^{k-2}  ) 
\\
= & \left(n + 3^{k-2} - 1 -  \left\lfloor\frac{n-1}{3}\right\rfloor \right)m    + 3^{k-2}\left( - 3n + 4 + 2\cdot 3^{k-2} + 4\left\lfloor\frac{n-1}{3}\right\rfloor \right)
\\
= &  B(n,m).
\end{align*}

\item $2\cdot 3^{k-1} \leq m \leq 2\cdot 3^{k-1} + 3^{k-2} .$

\begin{align*}
A(n,m) = &  S(3^{k-2}, 3^{k-2}, n - 2\cdot 3^{k-2},m) 
\\
= &  nm - 2A(m,3^{k-2})- A(m,n - 2\cdot 3^{k-2})
\\
= &  nm - 2\cdot 3^{2k-4} - (m - 3^{k-1})(n - 2\cdot 3^{k-2}) + 2m +  3^{k-3} \left(  m -42  + 2\cdot 3^{k-2} \right)
\\
& - \left\lfloor\frac{n -1}{3}\right\rfloor(14\cdot 3^{k-2} - 2m  )  +  2\cdot 3^{k-3}(14\cdot 3^{k-2} - 2m  ) 
\\
= &  \left(2 + 3^{k-2} + 2\left\lfloor\frac{n-1}{3}\right\rfloor  \right)m  + 3^{k-2}\left( 3n - 14 + 2\cdot 3^{k-2} - 14\left\lfloor\frac{n-1}{3}\right\rfloor    \right)
\\
= &  B(n,m).
\end{align*}

\item $ 2\cdot 3^{k-1} + 3^{k-2} \leq m \leq 2\cdot 3^{k-1} + 2\cdot 3^{k-2} .$

\begin{align*}
A(n,m) = &  S(3^{k-2}, 3^{k-2}, n - 2\cdot 3^{k-2},m) 
\\
= &  nm - 2A(m,3^{k-2})- A(m,n - 2\cdot 3^{k-2})
\\
= &  nm - 2\cdot 3^{2k-4} - (m - 3^{k-1})(n - 2\cdot 3^{k-2}) + m -  3^{k-2} \left( 7 -  3^{k-1} \right)
\\
& - \left\lfloor\frac{n -1}{3}\right\rfloor(14\cdot 3^{k-2} - 2m  )  +  2\cdot 3^{k-3}(14\cdot 3^{k-2} - 2m  ) 
\\
= &   \left(1 + 2\left\lfloor\frac{n-1}{2}\right\rfloor\right)m + 3^{k-2} \left( 3n - 7 + 3^{k} - 14\left\lfloor\frac{n-1}{2}\right\rfloor \right)
\\
= &  B(n,m).
\end{align*}

\item $  2\cdot 3^{k-1} + 2\cdot 3^{k-2} \leq m \leq 3^{k}.$

\begin{align*}
A(n,m) = &  S(3^{k-2}, 3^{k-2}, n - 2\cdot 3^{k-2},m) 
\\
= &  nm - 2A(m,3^{k-2})- A(m,n - 2\cdot 3^{k-2})
\\
= &  nm - 2\cdot 3^{2k-4} - 3^{k-1}(n - 2\cdot 3^{k-2}) - 2\cdot 3^{k-2}\left\lfloor\frac{n - 2\cdot 3^{k-2}}{2}\right\rfloor -  3^{k-2} + 3^{2k-3}
\\
= &  nm + 2\cdot 3^{k-2}\left\lfloor\frac{n-1}{2}\right\rfloor - 3^{k-2}( 5n -1 -  3^{k} )
\\
= &  B(n,m).
\end{align*}

\end{enumerate}
\end{proof}

\begin{lemma}\label{inductionlemma3}
Let $k \geq 4$ and suppose that $A(n',m') = B(n',m')$for all $n' \leq 3^{k-1}.$ Then

\begin{enumerate}

\item $ 3^{k-1} + 3^{k-2} \leq m \leq 2 \cdot 3^{k-1},$  $1+ 3^{k-1} \leq  x < 3^{k-1} + 3^{k-2}.$

$$A(m,x) = x\cdot 3^{k-1}+ \left\lfloor\frac{x-1}{3}\right\rfloor( m -4\cdot 3^{k-2}   )  +  m  - 3^{k-3}\left( m  + 12 + 2\cdot 3^{k-2} \right).$$

\item $  2 \cdot 3^{k-1} \leq m \leq  2 \cdot 3^{k-1} + 3^{k-2},$  $1+3^{k-1} \leq  x < 3^{k-1} + 3^{k-2}.$

$$A(m,x) = (m -  3^{k-1}  )x  + \left\lfloor\frac{x-1}{3}\right\rfloor (14\cdot 3^{k-2} - 2m )  - 2m + 3^{k-3} \left( - m + 42 - 2\cdot 3^{k-2}  \right).  $$

\item $  2 \cdot 3^{k-1} + 3^{k-2} \leq m \leq 2 \cdot 3^{k-1} + 2\cdot 3^{k-2} ,$  $1 + 3^{k-1} \leq  x < 3^{k-1} + 3^{k-2}.$

$$A(m,x) = (m - 3^{k-1}  )x + 3^{k-1}( 2\cdot 3^{k-1} - m  ). $$

\item $ 2 \cdot 3^{k-1} + 2\cdot 3^{k-2} \leq m \leq 3^{k},$  $1+3^{k-1} \leq  x < 3^{k-1} + 3^{k-2}.$

$$A(m,x) = (m - 3^{k-1}  )x + 3^{k-1}( 2\cdot 3^{k-1} - m  ). $$

\end{enumerate}

\end{lemma}

\begin{proof}
Follows directly from \cref{maintheorem2}.
\end{proof}

We next handle the case $n = 3^{k-1} + 2\cdot 3^{k-2} + 1.$

$$
B(3^{k-1} + 2\cdot 3^{k-2} + 1,m) = 
\begin{cases}
 13m\cdot 3^{k-3}   + 3^{k-2}\left( 1  - 19\cdot 3^{k-3}  \right), &\text{if $3^{k-1} + 3^{k-2} \leq  m \leq 2\cdot 3^{k-1}$} \\

\left(2 +  13\cdot 3^{k-3}   \right)m  - 3^{k-2}\left( 11 - 19\cdot 3^{k-2}   \right) , &\text{if $2\cdot 3^{k-1} \leq m \leq 2\cdot 3^{k-1} + 3^{k-2}$} \\

5m\cdot 3^{k-2} + 3^{k-1}  - 11\cdot 3^{2k-4} , &\text{if  $ 2\cdot 3^{k-1} + 3^{k-2} \leq m \leq 2\cdot 3^{k-1} + 2\cdot 3^{k-2} $}\\

5m\cdot 3^{k-2} + m  - 11\cdot 3^{2k-4}  - 4\cdot  3^{k-2} , &\text{if $  2\cdot 3^{k-1} + 2\cdot 3^{k-2} \leq m $ } 

\end{cases}
$$

\begin{prop}\label{maintheorem3b}
Let $k \geq 4$ and suppose that $A(n',m') = B(n',m')$for all $n' \leq 3^{k-1}.$ Let $n = 3^{k-1} + 2\cdot 3^{k-2} + 1$ and $3^{k-1} + 3^{k-2} \leq m \leq 3^{k}. $ Then $(3^{k-2}, 3^{k-2} + 2, 3^{k-1} - 1)$ is an optimal move for $(n,m)$ and $A(n,m) = B(n,m).$

\end{prop}

\begin{proof}

Let $(a,b,c)$ be an optimal move. As in \cref{optimalform3} we may assume that $3^{k-2} \leq a \leq b\leq c.$ 

Suppose $a \geq 3^{k-2} + \alpha$ where $\alpha>0.$ Then $c \leq 3^{k-2} - 2\alpha + 1.$ Thus we may assume that $a= 3^{k-2}$ and $c\leq 3^{k-2}$ by the same reasoning as in \cref{optimalform3}.

By switching to $(a,b-3,c+3)$ if needed we may assume that $3^{k-1} - 2 \leq c \leq 3^{k-1} + 1.$

We break up the rest of the proof into cases:

\begin{enumerate}

\item $ 3^{k-1} + 3^{k-2} \leq m \leq   2\cdot 3^{k-1}.$  

By \cref{inductionlemma2,inductionlemma3} we have

$$
A(m,c) - A(m,c-1) = 
\begin{cases}
m - 3^{k-2}, &\text{if $c = 3^{k-1} + 1$ } \\

3^{k-1}, &\text{if $c = 3^{k-1} $} \\

3^{k-1}, &\text{if $c = 3^{k-1} - 1$}

\end{cases}
$$

$$
A(m,b+1) - A(m,b) = 
\begin{cases}
3^{k-1}, &\text{if $b = 3^{k-2} + 2$ } \\

3^{k-1}, &\text{if $b = 3^{k-2} + 1 $} \\

m - 3^{k-2}, &\text{if $b = 3^{k-2} $}

\end{cases}
$$

Thus $A(m,c) - A(m,c-1) = A(m,b+1) - A(m,b)$ for all $b,c$ in our range and we may switch to $(3^{k-2}, 3^{k-2} + 2, 3^{k-1} - 1).$

Then

\begin{align*}
A(n,m) = &  S(3^{k-2}, 3^{k-2} + 2, 3^{k-1} - 1,m) 
\\
= &  ( 5\cdot 3^{k-2} + 1)m - A(m,3^{k-2})- A(m,3^{k-2} + 2) - A(m,3^{k-1} - 1) 
\\
= &  5m\cdot 3^{k-2} + m - 3^{2k - 4} - 3^{2k - 4} - m - 2\cdot 3^{k-2} - 13\cdot 3^{2k-5} + 3^{k-1}  - 2m\cdot 3^{k-3}
\\
= &  13m\cdot 3^{k-3}   + 3^{k-2}\left( 1  - 19\cdot 3^{k-3}  \right)
\\
= &  B(n,m).
\end{align*}

\item $  2\cdot 3^{k-1} \leq m \leq 2\cdot 3^{k-1} + 3^{k-2}.$

By \cref{inductionlemma2,inductionlemma3} we have

$$
A(m,c) - A(m,c-1) = 
\begin{cases}
11\cdot 3^{k-2} - m, &\text{if $c = 3^{k-1} + 1$ } \\

m - 3^{k-1}, &\text{if $c = 3^{k-1} $} \\

m - 3^{k-1}, &\text{if $c = 3^{k-1} - 1$}

\end{cases}
$$

$$
A(m,b+1) - A(m,b) = 
\begin{cases}
m - 3^{k-1}, &\text{if $b = 3^{k-2} + 2$ } \\

3^{k-1}, &\text{if $b = 3^{k-2} + 1 $} \\

5\cdot 3^{k-2}, &\text{if $b = 3^{k-2} $}

\end{cases}
$$

If $c = 3^{k-1}$ and $b = 3^{k-2} + 1.$ Then $A(m,c+1)-A(m,c) = 11\cdot 3^{k-2} - m \leq 5\cdot 3^{k-2} = A(m,b) - A(m,b-1). $

If $c = 3^{k-1} + 1$ and $b = 3^{k-2} .$ Then $A(m,c)-A(m,c-2) = 8\cdot 3^{k-2} = A(m,b) - A(m,b-1). $

Finally suppose $c = 3^{k-1} - 2$ and $b = 3^{k-2} + 3 .$ Then $A(m,c+1)-A(m,c) = m - 3^{k-1} = A(m,b) - A(m,b-1). $

Thus we may switch to $(3^{k-2}, 3^{k-2} + 2, 3^{k-1} - 1).$

Then

\begin{align*}
A(n,m) = &  S(3^{k-2}, 3^{k-2} + 2, 3^{k-1} - 1,m) 
\\
= &  ( 5\cdot 3^{k-2} + 1)m - A(m,3^{k-2})- A(m,3^{k-2} + 2) - A(m,3^{k-1} - 1) 
\\
= &  5m\cdot 3^{k-2} + m - 3^{2k - 4} -  8\cdot 3^{k-2} - 3^{2k - 4} - 2m\cdot 3^{k-3} + m - 3^{k-1} - 13\cdot 3^{2k-5} 
\\
= &  13m\cdot 3^{k-3} + 2m  -  11\cdot 3^{k-2} - 19\cdot 3^{2k-5} 
\\
= &  B(n,m).
\end{align*}

\item $ 2\cdot 3^{k-1} + 3^{k-2} \leq m \leq 2\cdot 3^{k-1} + 2\cdot 3^{k-2}.$

By \cref{inductionlemma2,inductionlemma3} we have

$$
A(m,c) - A(m,c-1) = 
\begin{cases}
m - 3^{k-1}, &\text{if $c = 3^{k-1} + 1$ } \\

11\cdot 3^{k-2} - m, &\text{if $c = 3^{k-1} $} \\

m - 3^{k-1}, &\text{if $c = 3^{k-1} - 1$}

\end{cases}
$$

$$
A(m,b+1) - A(m,b) = 
\begin{cases}
m - 3^{k-1}, &\text{if $b = 3^{k-2} + 2$ } \\

3^{k-1}, &\text{if $b = 3^{k-2} + 1 $} \\

5\cdot 3^{k-2}, &\text{if $b = 3^{k-2} $}

\end{cases}
$$

Suppose $c = 3^{k-1}$ and $b = 3^{k-2} + 1.$ Then $A(m, c+1) - A(m,c) = m - 3^{k-1} \leq 5\cdot 3^{k-2} = A(m,b) -A(m,b-1). $

Suppose $c = 3^{k-1}+1$ and $b = 3^{k-2}.$ Then $A(m, c) - A(m,c-2) = 8\cdot 3^{k-2} = A(m,b+2) -A(m,b). $

Finally suppose $c = 3^{k-1}-2$ and $b = 3^{k-2} + 3.$ Then $A(m, c + 1) - A(m,c) = m - 3^{k-1} = A(m,b) -A(m,b-1). $

Thus we may switch to $(3^{k-2}, 3^{k-2} + 2, 3^{k-1} - 1).$

Then

\begin{align*}
A(n,m) = &  S(3^{k-2}, 3^{k-2} + 2, 3^{k-1} - 1,m) 
\\
= &  ( 5\cdot 3^{k-2} + 1)m - A(m,3^{k-2})- A(m,3^{k-2} + 2) - A(m,3^{k-1} - 1) 
\\
= &  5m\cdot 3^{k-2} + m - 3^{2k - 4} - 8\cdot 3^{k-2} - 3^{2k - 4} - m  + 11\cdot 3^{k-2} - 3^{2k-2} 
\\
= &  5m\cdot 3^{k-2} + 3^{k-1} - 11\cdot 3^{2k - 4} 
\\
= &  B(n,m).
\end{align*}

\item $2\cdot 3^{k-1} + 2\cdot 3^{k-2} \leq m \leq 3^{k}.$

By \cref{inductionlemma2,inductionlemma3} we have

$$
A(m,c) - A(m,c-1) = 
\begin{cases}
m - 3^{k-1}, &\text{if $c = 3^{k-1} + 1$ } \\

3^{k-1}, &\text{if $c = 3^{k-1} $} \\

5\cdot 3^{k-2}, &\text{if $c = 3^{k-1} - 1$}

\end{cases}
$$

$$
A(m,b+1) - A(m,b) = 
\begin{cases}
5\cdot 3^{k-2}, &\text{if $b = 3^{k-2} + 2$ } \\

3^{k-1}, &\text{if $b = 3^{k-2} + 1 $} \\

5\cdot 3^{k-2}, &\text{if $b = 3^{k-2} $}

\end{cases}
$$

Suppose $c = 3^{k-1} + 1$ and $b = 3^{k-2}.$ Then $A(m, c) - A(m,c-1) \geq 5\cdot 3^{k-2} = A(m,b+1) -A(m,b). $

Suppose $c = 3^{k-1} $ and $b = 3^{k-2}+1.$ Then $A(m, c) - A(m,c-2) = 8\cdot 3^{k-2} = A(m,b+2) -A(m,b). $

Suppose $c = 3^{k-1} -2 $ and $b = 3^{k-2}+3.$ Then $A(m, c+1) - A(m,c) = 5\cdot 3^{k-2} = A(m,b) -A(m,b-1). $

Thus we may switch to $(3^{k-2}, 3^{k-2} + 2, 3^{k-1} - 1).$

Then

\begin{align*}
A(n,m) = &  S(3^{k-2}, 3^{k-2} + 2, 3^{k-1} - 1,m) 
\\
= &  ( 5\cdot 3^{k-2} + 1)m - A(m,3^{k-2})- A(m,3^{k-2} + 2) - A(m,3^{k-1} - 1) 
\\
= &  5m\cdot 3^{k-2} + m - 3^{2k - 4} - 8\cdot 3^{k-2} - 3^{2k - 4} -  3^{2k-2} + 4\cdot 3^{k-2}
\\
= &  5m\cdot 3^{k-2} + m  - 4\cdot 3^{k-2} - 11\cdot 3^{2k - 4}  
\\
= &  B(n,m).
\end{align*}

\end{enumerate}
    
\end{proof}

\begin{prop}\label{maintheorem3c}
Let $k \geq 4$ and suppose that $A(n',m') = B(n',m')$for all $n' \leq 3^{k-1}.$ Let $3^{k-1} + 2\cdot 3^{k-2} + 1 < n \leq 2\cdot 3^{k-1} + 3^{k-2}$ and $ 3^{k-1} + 3^{k-2} \leq m  \leq 2\cdot 3^{k-1}  + 3^{k-2}. $ Then $(3^{k-2}, n - 3^{k-2} - 3^{k-1}, 3^{k-1} )$ is an optimal move for $(n,m)$ and $A(n,m) = B(n,m).$

\end{prop}

\begin{proof}

Let $(a,b,c)$ be an optimal move. As in \cref{optimalform3} we may assume that $3^{k-2} \leq a \leq b\leq c.$ 

We split the proof in to two cases

\begin{enumerate}
\item  $ 3^{k-1} + 3^{k-2} \leq m \leq   2\cdot 3^{k-1}. $

For $1+ 3^{k-1} \leq  x < 3^{k-1} + 3^{k-2}$ we have the following by \cref{inductionlemma3}

$$A(m,x) = x\cdot 3^{k-1}+ \left\lfloor\frac{x-1}{3}\right\rfloor( m -4\cdot 3^{k-2}   )  +  m  - 3^{k-3}\left( m  + 12 + 2\cdot 3^{k-2} \right).$$

For $3^{k-1} + 3^{k-2} \leq  x \leq 3^{k-1} + 2\cdot 3^{k-2}$ we have the following consequence of \cref{maintheorem3}

$$A(m,x) = \left(m + 3^{k-2} - 1 - \left\lfloor\frac{m-1}{3}\right\rfloor \right)x    + 3^{k-2}\left( - 3m + 4 + 2\cdot 3^{k-2} + 4\left\lfloor\frac{m-1}{3}\right\rfloor \right).$$

Recall that for $1+3^{k-2}\leq x \leq 3^{k-1}$ we have

$$A(m,x) = x\cdot 3^{k-1} +  \left\lfloor\frac{x-1}{3}\right\rfloor(m - 4\cdot 3^{k-2}  ) + m - 3^{k-3} \left(  m + 12  + 2\cdot 3^{k-2} \right). $$

If $c \geq 3^{k-1} + 3 $ then $b\leq 3^{k-1} - 3$ and $A(m,b+3)-A(m,b) \leq A(m,c) - A(m,c-3).$ Note that $A(m, 1 + 3^{k-1}) - A(m, 3^{k}) = m - 3^{k-2} \geq 3^{k-1}$. Therefore if $c = 2 + 3^{k-1}$ we have  $A(m,b+2)-A(m,b) \leq A(m,c) - A(m,c-2)$ and if $c = 1 + 3^{k-1}$ we have  $A(m,b+1)-A(m,b) \leq A(m,c) - A(m,c-1).$ Thus we may assume that $c\leq 3^{k-1}.$

Using a similar argument to above we may assume that $c\geq 3^{k-1} - 2.$ Since $A(m, 3^{k-1}) - A(m, 3^{k} - 1) = A(m, 3^{k-1}-1) - A(m, 3^{k} - 2) = 3^{k-1} $ we may raise $c$ and lower $a$ or $b.$ Thus we can assume that $c = 3^{k-1}.$ 

In a similar way to how we lowered $c$ and raised $b$ above we now lower $a$ and raise $b$ until $a = 3^{k-2}.$ Thus arriving at the desired optimal move.

Thus we have 
\begin{align*}
A(n,m) = &  S(3^{k-2}, n - 4\cdot 3^{k-2}, 3^{k-1},m) 
\\
= &  nm - A(m,3^{k-2})  - A(m,3^{k-1}) - A(m,n - 4\cdot 3^{k-2})
\\
= &  nm -  3^{2k-4} - 2m\cdot 3^{k-3} - 13\cdot 3^{2k-5} -  \left\lfloor\frac{n -1}{3}\right\rfloor(m - 4\cdot 3^{k-2}  )
\\
& -   3^{k-1}(n - 4\cdot 3^{k-2})  - m + 3^{k-3} \left(  m + 12  + 2\cdot 3^{k-2} \right) + 4\cdot 3^{k-3}(m - 4\cdot 3^{k-2}  )
\\
= & \left(n + 3^{k-2} - 1 -  \left\lfloor\frac{n-1}{3}\right\rfloor \right)m    + 3^{k-2}\left( - 3n + 4 + 2\cdot 3^{k-2} + 4\left\lfloor\frac{n-1}{3}\right\rfloor \right)
\\
= &  B(n,m).
\end{align*}

\item $ 2\cdot 3^{k-1}\leq m \leq   2\cdot 3^{k-1}  + 3^{k-2} . $

For $1+ 3^{k-1} \leq  x < 3^{k-1} + 3^{k-2}$ we have the following by \cref{inductionlemma3}

$$A(m,x) = (m -  3^{k-1}  )x  + \left\lfloor\frac{x-1}{3}\right\rfloor (14\cdot 3^{k-2} - 2m )  - 2m + 3^{k-3} \left( - m + 42 - 2\cdot 3^{k-2}  \right).  $$

For $3^{k-1} + 3^{k-2} \leq  x \leq 3^{k-1} + 2\cdot 3^{k-2}$ we have the following consequence of part (1).

$$A(m,x) = \left(m + 3^{k-2} - 1 - \left\lfloor\frac{m-1}{3}\right\rfloor \right)x    + 3^{k-2}\left( - 3m + 4 + 2\cdot 3^{k-2} + 4\left\lfloor\frac{m-1}{3}\right\rfloor \right).$$

Recall that for $1+3^{k-2}\leq x \leq 3^{k-1}$ we have

$$A(m,x) = (m - 3^{k-1} )x +\left\lfloor\frac{x-1}{3}\right\rfloor(14\cdot 3^{k-2}  -2m) - 2m - 3^{k-3} \left( m - 42  + 2\cdot 3^{k-2} \right). $$

We can obtain the desired optimal move by following the proof of part $(1)$.

Thus we have 
\begin{align*}
A(n,m) = &  S(3^{k-2}, n - 4\cdot 3^{k-2}, 3^{k-1},m) 
\\
= &  nm - A(m,3^{k-2})  - A(m,3^{k-1}) - A(m,n - 4\cdot 3^{k-2})
\\
= &  nm -  3^{2k-4} - 2m\cdot 3^{k-3} - 13\cdot 3^{2k-5} - 4\cdot 3^{k-3}(14\cdot 3^{k-2}  -2m)
\\
& -   (m - 3^{k-1} )(n - 4\cdot 3^{k-2}) +\left\lfloor\frac{n -1}{3}\right\rfloor(14\cdot 3^{k-2}  -2m) - 2m - 3^{k-3} \left( m - 42  + 2\cdot 3^{k-2} \right)
\\
= & \left(2 + 3^{k-2} + 2\left\lfloor\frac{n-1}{3}\right\rfloor  \right)m  + 3^{k-2}\left( 3n - 14 + 2\cdot 3^{k-2} - 14\left\lfloor\frac{n-1}{3}\right\rfloor    \right)
\\
= &  B(n,m).
\end{align*}

\end{enumerate}
\end{proof}

We now partially handle the case $n = 2\cdot 3^{k-1} +  3^{k-2} + 1.$

$$
B(2\cdot 3^{k-1} +  3^{k-2} + 1,m) = 
\begin{cases}
17m\cdot 3^{k-3}   + 3^{k-2}  - 29\cdot 3^{2k-5} , &\text{if $3^{k-1} + 3^{k-2} \leq  m \leq 2\cdot 3^{k-1}$} \\

\left(2 + 17\cdot 3^{k-3} \right)m  - 11\cdot 3^{k-2}  - 29\cdot 3^{2k-5}, &\text{if $2\cdot 3^{k-1} \leq m \leq 2\cdot 3^{k-1} + 3^{k-2}$} \\

7m\cdot 3^{k-2} + 11\cdot 3^{k-2} - 19\cdot 3^{2k-4} , &\text{if  $ 2\cdot 3^{k-1} + 3^{k-2} \leq m \leq 2\cdot 3^{k-1} + 2\cdot 3^{k-2} $}\\

 7m\cdot 3^{k-2} + m  = 5\cdot 3^{k-2} -  19\cdot 3^{2k-4} , &\text{if $  2\cdot 3^{k-1} + 2\cdot 3^{k-2} \leq m $ } 

\end{cases}
$$

\begin{prop}\label{maintheorem3d}
Let $k \geq 4$ and suppose that $A(n',m') = B(n',m')$for all $n' \leq 3^{k-1}.$ Let $n = 2\cdot 3^{k-1} + 3^{k-2} + 1$ and $3^{k-1} + 3^{k-2} \leq m \leq 2\cdot 3^{k-1} + 3^{k-2}. $ Then $A(n,m) = B(n,m).$

\end{prop}

\begin{proof}

Following the proof of \cref{maintheorem3b} we can show that  $(3^{k-2} + 2,  3^{k-1} - 1, 3^{k-1} )$ is optimal.

We split up the rest of the proof into two cases:

\begin{enumerate}

\item $ 3^{k-1} + 3^{k-2} \leq m \leq   2\cdot 3^{k-1}.$

\begin{align*}
A(n,m) = &  S(3^{k-2} + 2,  3^{k-1} - 1, 3^{k-1},m) 
\\
= &  ( 7\cdot 3^{k-2} + 1)m - A(m,3^{k-1})- A(m,3^{k-2} + 2) - A(m,3^{k-1} - 1) 
\\
= &  7m\cdot 3^{k-2} + m - 2m\cdot 3^{k-3} - 13\cdot 3^{2k-5} -  3^{2k - 4} - m - 2\cdot 3^{k-2} - 13\cdot 3^{2k-5} + 3^{k-1}  - 2m\cdot 3^{k-3}
\\
= &  17m\cdot 3^{k-3}   + 3^{k-2}  - 29\cdot 3^{2k-5} 
\\
= &  B(n,m).
\end{align*}

\item $  2\cdot 3^{k-1} \leq m \leq 2\cdot 3^{k-1} + 3^{k-2}.$

\begin{align*}
A(n,m) = &  S(3^{k-2} + 2,  3^{k-1} - 1, 3^{k-1},m) 
\\
= &  ( 7\cdot 3^{k-2} + 1)m - A(m,3^{k-1})- A(m,3^{k-2} + 2) - A(m,3^{k-1} - 1) 
\\
= &  7m\cdot 3^{k-2} + m -  2m\cdot 3^{k-3} - 13\cdot 3^{2k-5} - 8\cdot 3^{k-2} -  3^{2k - 4} -  2m\cdot 3^{k-3} + m - 3^{k-1} - 13\cdot 3^{2k-5} 
\\
= & \left(2 + 17\cdot 3^{k-3} \right)m  - 11\cdot 3^{k-2}  - 29\cdot 3^{2k-5}
\\
= &  B(n,m).
\end{align*}

\end{enumerate}
    
\end{proof}

\begin{prop}\label{maintheorem3e}
Let $k \geq 4$ and suppose that $A(n',m') = B(n',m')$for all $n' \leq 3^{k-1}.$ Let $ 2\cdot  3^{k-1} + 3^{k-2} + 1 < n \leq 3^k$ and $ 3^{k-1} + 3^{k-2} \leq m  \leq 2\cdot 3^{k-1}  + 3^{k-2}. $ Then $(3^{k-2}, n - 3^{k-2} - 3^{k-1}, 3^{k-1} )$ is an optimal move for $(n,m)$ and $A(n,m) = B(n,m).$

\end{prop}

\begin{proof}

Let $(a,b,c)$ be an optimal move. As in \cref{optimalform3} we may assume that $3^{k-2} \leq a \leq b\leq c.$ 

We split the proof in to two cases

\begin{enumerate}
\item  $ 3^{k-1} + 3^{k-2} \leq m \leq   2\cdot 3^{k-1}. $

For $1+ 3^{k-1} \leq  x < 3^{k-1} + 3^{k-2}$ we have the following by \cref{inductionlemma3}

$$A(m,x) = x\cdot 3^{k-1}+ \left\lfloor\frac{x-1}{3}\right\rfloor( m -4\cdot 3^{k-2}   )  +  m  - 3^{k-3}\left( m  + 12 + 2\cdot 3^{k-2} \right).$$

For $3^{k-1} + 3^{k-2} \leq  x \leq 3^{k-1} + 2\cdot 3^{k-2}$ we have the following consequence of \cref{maintheorem3}

$$A(m,x) = \left(m + 3^{k-2} - 1 - \left\lfloor\frac{m-1}{3}\right\rfloor \right)x    + 3^{k-2}\left( - 3m + 4 + 2\cdot 3^{k-2} + 4\left\lfloor\frac{m-1}{3}\right\rfloor \right).$$

For $ 3^{k-1} + 2\cdot 3^{k-2} \leq  x \leq 2\cdot 3^{k-1}  + 3^{k-2}$ we have the following consequence of \cref{maintheorem3c}

$$A(m,x) = \left(2 + 3^{k-2} + 2\left\lfloor\frac{m-1}{3}\right\rfloor  \right)x  + 3^{k-2}\left( 3m - 14 + 2\cdot 3^{k-2} - 14\left\lfloor\frac{m-1}{3}\right\rfloor    \right) $$

Recall that for $1+3^{k-2}\leq x \leq 3^{k-1}$ we have

$$A(m,x) = x\cdot 3^{k-1} +  \left\lfloor\frac{x-1}{3}\right\rfloor(m - 4\cdot 3^{k-2}  ) + m - 3^{k-3} \left(  m + 12  + 2\cdot 3^{k-2} \right). $$

Since $n\leq 3^{k}$ and $3^{k-2}\leq a \leq b$ we have $c\leq 2\cdot 3^{k-1}  + 3^{k-2}.$

If $c \geq 3^{k-1} + 3 $ then $b\leq 3^{k-1} - 3.$ Then by \cref{allgaps} and \cref{maintheorem3c} we have $A(m,b+3)-A(m,b) \leq A(m,c) - A(m,c-3).$ Also from  \cref{maintheorem3c} we have $A(m, 1 + 3^{k-1}) - A(m, 3^{k}) = m - 3^{k-2} \geq 3^{k-1}$. Continuing the argument from the proof of \cref{maintheorem3c} we may assume that $c = 3^{k-1}.$ Swapping $c$ for $b$ and $b$ for $a$ and repeating this argument shows we may assume that $b = 3^{k-1}.$ Thus arriving at the desired optimal move.

Thus we have 
\begin{align*}
A(n,m) = &  S( 3^{k-1}, 3^{k-1}, n - 2\cdot 3^{k-1},  m) 
\\
= &  nm -  2\cdot A(m,3^{k-1}) - A(m,  n - 2\cdot 3^{k-1})
\\
= &  nm - 4m\cdot 3^{k-3} - 26\cdot 3^{2k-5} - 3^{k-1}(n - 2\cdot 3^{k-1}) +  2\cdot 3^{k-2}(m - 4\cdot 3^{k-2}  )
\\
& -  \left\lfloor\frac{n -1}{3}\right\rfloor(m - 4\cdot 3^{k-2}  ) - m + 3^{k-3} \left(  m + 12  + 2\cdot 3^{k-2} \right) 
\\
= & \left(n + 3^{k-2} - 1 -  \left\lfloor\frac{n-1}{3}\right\rfloor \right)m    + 3^{k-2}\left( - 3n + 4 + 2\cdot 3^{k-2} + 4\left\lfloor\frac{n-1}{3}\right\rfloor \right)
\\
= &  B(n,m).
\end{align*}

\item $ 2\cdot 3^{k-1}\leq m \leq   2\cdot 3^{k-1}  + 3^{k-2} . $

We can obtain the desired optimal move by following the proof of part $(1)$.

Thus we have 
\begin{align*}
A(n,m) = &  S( 3^{k-1}, 3^{k-1}, n - 2\cdot 3^{k-1},  m) 
\\
= &  nm -  2\cdot A(m,3^{k-1}) - A(m,  n - 2\cdot 3^{k-1})
\\
= &  nm -  4m\cdot 3^{k-3} - 26\cdot 3^{2k-5}  -  (m - 3^{k-1} )(n - 2\cdot 3^{k-1}) +   2\cdot 3^{k-2} (14\cdot 3^{k-2}  -2m)
\\
&  -\left\lfloor\frac{n -1}{3}\right\rfloor(14\cdot 3^{k-2}  -2m) + 2m + 3^{k-3} \left( m - 42  + 2\cdot 3^{k-2} \right)
\\
= & \left(2 + 3^{k-2} + 2\left\lfloor\frac{n-1}{3}\right\rfloor  \right)m  + 3^{k-2}\left( 3n - 14 + 2\cdot 3^{k-2} - 14\left\lfloor\frac{n-1}{3}\right\rfloor    \right)
\\
= &  B(n,m).
\end{align*}

\end{enumerate}
\end{proof}

\begin{prop}\label{maintheorem3f}
Let $k \geq 4$ and suppose that $A(n',m') = B(n',m')$for all $n' \leq 3^{k-1}.$ Let $  3^{k-1} + 2\cdot 3^{k-2} + 1 < n \leq 3^k$ and $ 2\cdot 3^{k-1}  + 3^{k-2} \leq m \leq 2\cdot 3^{k-1} + 2\cdot 3^{k-2}. $ Then $A(n,m) = B(n,m).$

\end{prop}

\begin{proof}

\cref{maintheorem3e} gives us the last formulas needed to handle this case. Namely for $3^{k-1} + 3^{k-2} \leq  x \leq 3^{k-1} + 2\cdot 3^{k-2}$ we have the following consequence of \cref{maintheorem3}

$$A(m,x) = \left(m + 3^{k-2} - 1 - \left\lfloor\frac{m-1}{3}\right\rfloor \right)x    + 3^{k-2}\left( - 3m + 4 + 2\cdot 3^{k-2} + 4\left\lfloor\frac{m-1}{3}\right\rfloor \right),$$

and for $ 3^{k-1} + 2\cdot 3^{k-2} \leq  x \leq 2\cdot 3^{k-1}  + 3^{k-2}$ we have

$$A(m,x) = \left(2 + 3^{k-2} + 2\left\lfloor\frac{m-1}{3}\right\rfloor  \right)x  + 3^{k-2}\left( 3m - 14 + 2\cdot 3^{k-2} - 14\left\lfloor\frac{m-1}{3}\right\rfloor    \right). $$

Following the proofs of \cref{maintheorem3c,maintheorem3d,maintheorem3e} shows that the optimal moves described there are also optimal in this case.

We split the remainder of the proof in to cases

\begin{enumerate}
\item  $3^{k-1} + 2\cdot 3^{k-2} + 1 < n \leq 2\cdot  3^{k-1} + 3^{k-2} $

\begin{align*}
A(n,m) = &  S(3^{k-2}, n - 4\cdot 3^{k-2}, 3^{k-1},m) 
\\
= &  nm - A(m,3^{k-2})  - A(m,3^{k-1}) - A(m,n - 4\cdot 3^{k-2})
\\
= &  nm - 2\cdot 3^{2k-4} -  3^{2k-2} -  (m - 3^{k-1} )(n - 4\cdot 3^{k-2}) +   4\cdot 3^{k-3}( 14\cdot 3^{k-2} - 2m  )
\\
&  - \left\lfloor\frac{n -1}{2}\right\rfloor( 14\cdot 3^{k-2} - 2m  ) + m - 3^{k-2} \left( 7 - 3^{k-1}\right)
\\
= &   \left(1 + 2\left\lfloor\frac{n-1}{2}\right\rfloor\right)m + 3^{k-2} \left( 3n - 7 + 3^{k} - 14\left\lfloor\frac{n-1}{2}\right\rfloor \right)
\\
= &  B(n,m).
\end{align*}

\item $n = 2\cdot  3^{k-1} + 3^{k-2} + 1$

\begin{align*}
A(n,m) = &  S(3^{k-2} + 2,  3^{k-1} - 1, 3^{k-1},m) 
\\
= &  ( 7\cdot 3^{k-2} + 1)m - A(m,3^{k-1}) - A(m,3^{k-1} - 1) - A(m,3^{k-2} + 2)
\\
= &  7m\cdot 3^{k-2} + m - 3^{2k-2} -  m  + 11\cdot 3^{k-2} + 3^{2k-2} -  8\cdot 3^{k-2} -  3^{2k - 4} 
\\
= & 7m\cdot 3^{k-2} + 11\cdot 3^{k-2} - 19\cdot 3^{2k-4} 
\\
= &  B(n,m).
\end{align*}

\item $2\cdot  3^{k-1} + 3^{k-2} + 1 < n \leq 3^k$

\begin{align*}
A(n,m) = &  S( 3^{k-1}, 3^{k-1}, n - 2\cdot 3^{k-1},  m) 
\\
= &  nm -  2\cdot A(m,3^{k-1}) - A(m,  n - 2\cdot 3^{k-1})
\\
= &  nm - 2\cdot 3^{2k-2} - (m - 3^{k-1} )( n - 2\cdot 3^{k-1}) + 3^{k-1}( 14\cdot 3^{k-2} - 2m  )
\\
&   - \left\lfloor\frac{ n -1}{2}\right\rfloor( 14\cdot 3^{k-2} - 2m  ) + m - 3^{k-2} \left( 7 - 3^{k-1}\right) 
\\
= &  \left(1 + 2\left\lfloor\frac{n-1}{2}\right\rfloor\right)m + 3^{k-2} \left( 3n - 7 + 3^{k} - 14\left\lfloor\frac{n-1}{2}\right\rfloor \right)
\\
= &  B(n,m).
\end{align*}

\end{enumerate}
    
\end{proof}

\begin{prop}
Let $k \geq 4$ and suppose that $A(n',m') = B(n',m')$for all $n' \leq 3^{k-1}.$ Let $  3^{k-1} + 2\cdot 3^{k-2} + 1 < n \leq 3^k$ and $ 2\cdot 3^{k-1}  + 2\cdot 3^{k-2} \leq m \leq  3^{k}. $ Then $A(n,m) = B(n,m).$

\end{prop}

\begin{proof}

\cref{maintheorem3e} gives us the last formulas needed to handle this case. Namely for $3^{k-1} + 3^{k-2} \leq  x \leq 3^{k-1} + 2\cdot 3^{k-2}$ we have the following consequence of \cref{maintheorem3}

$$A(m,x) = \left(m + 3^{k-2} - 1 - \left\lfloor\frac{m-1}{3}\right\rfloor \right)x    + 3^{k-2}\left( - 3m + 4 + 2\cdot 3^{k-2} + 4\left\lfloor\frac{m-1}{3}\right\rfloor \right),$$

and for $ 3^{k-1} + 2\cdot 3^{k-2} \leq  x \leq 2\cdot 3^{k-1}  + 3^{k-2}$ we have

$$A(m,x) = \left(2 + 3^{k-2} + 2\left\lfloor\frac{m-1}{3}\right\rfloor  \right)x  + 3^{k-2}\left( 3m - 14 + 2\cdot 3^{k-2} - 14\left\lfloor\frac{m-1}{3}\right\rfloor    \right). $$

Following the proofs of Theorems \cref{maintheorem3c,maintheorem3d,maintheorem3e} shows that the optimal moves described there are also optimal in this case.

We split the remainder of the proof in to cases

\begin{enumerate}
\item  $3^{k-1} + 2\cdot 3^{k-2} + 1 < n \leq 2\cdot  3^{k-1} + 3^{k-2} $

\begin{align*}
A(n,m) = &  S(3^{k-2}, n - 4\cdot 3^{k-2}, 3^{k-1},m) 
\\
= &  nm - A(m,3^{k-2})  - A(m,3^{k-1}) - A(m,n - 4\cdot 3^{k-2})
\\
= &  nm - 2\cdot 3^{2k-4} -  3^{2k-2} + 3^{k-2} + 3^{k-1}(n - 4\cdot 3^{k-2}) 
\\
&  - 2\cdot 3^{k-2}\left\lfloor\frac{n - 4\cdot 3^{k-2}}{2}\right\rfloor - 3^{k-2} + 3^{2k-3}
\\
= &   nm + 2\cdot 3^{k-2}\left\lfloor\frac{n-1}{2}\right\rfloor - 3^{k-2}( 5n -1 -  3^{k} )
\\
= &  B(n,m).
\end{align*}

\item $n = 2\cdot  3^{k-1} + 3^{k-2} + 1$

\begin{align*}
A(n,m) = &  S(3^{k-2} + 2,  3^{k-1} - 1, 3^{k-1},m) 
\\
= &  ( 7\cdot 3^{k-2} + 1)m - A(m,3^{k-1}) - A(m,3^{k-1} - 1) - A(m,3^{k-2} + 2)
\\
= &  7m\cdot 3^{k-2} + m - 3^{2k-2} + 3^{k-2} -  3^{2k-2} + 4\cdot 3^{k-2} - 8\cdot 3^{k-2} - 3^{2k - 4}
\\
= & 7m\cdot 3^{k-2} + m  = 5\cdot 3^{k-2} -  19\cdot 3^{2k-4}
\\
= &  B(n,m).
\end{align*}

\item $2\cdot  3^{k-1} + 3^{k-2} + 1 < n \leq 3^k$

\begin{align*}
A(n,m) = &  S( 3^{k-1}, 3^{k-1}, n - 2\cdot 3^{k-1},  m) 
\\
= &  nm -  2\cdot A(m,3^{k-1}) - A(m,  n - 2\cdot 3^{k-1})
\\
= &  nm -   2\cdot 3^{2k-2} +  2\cdot 3^{k-2} +  3^{k-1} ( n - 2\cdot 3^{k-1})
\\
&   -  2\cdot 3^{k-2}\left\lfloor\frac{ n - 2\cdot 3^{k-1}}{2}\right\rfloor - 3^{k-2} + 3^{2k-3} 
\\
= &   nm + 2\cdot 3^{k-2}\left\lfloor\frac{n-1}{2}\right\rfloor - 3^{k-2}( 5n -1 -  3^{k} )
\\
= &  B(n,m).
\end{align*}

\end{enumerate}
    
\end{proof}

\begin{lemma}\label{modsumhelper}
Let $k\geq 4$ and $ a,b,c \in\mathbb{N}$ and set $n = a + b + c.$ Suppose that at least two of $a,b$ and $c$ are odd. Then

$$\left\lfloor\frac{n-1}{2}\right\rfloor = \left\lfloor\frac{a+1}{2}\right\rfloor + \left\lfloor\frac{b+1}{2}\right\rfloor + \left\lfloor\frac{c+1}{2}\right\rfloor - 2.$$

\end{lemma}

\begin{proof}

Without loss of generality we assume that $a = 1 \mod 2$ and $b = 1 \mod 2$. Then by using

$$(a+1) = 0 \mod 2, (b+1) = 0 \mod 2, -4 = 0 \mod 2,$$

we have

\begin{align*}
\left\lfloor\frac{n-1}{2}\right\rfloor & = \left\lfloor\frac{(a+1)+(b+1)+(c+1)-4}{2}\right\rfloor
\\
& =\left\lfloor\frac{a+1}{2}\right\rfloor + \left\lfloor\frac{b+1}{2}\right\rfloor + \left\lfloor\frac{c+1}{2}\right\rfloor - 2.
\end{align*}
\end{proof}

\begin{prop}
Let $k \geq 4$ and suppose that $A(n',m') = B(n',m')$for all $n' \leq 3^{k-1}.$ Let $  3^{k-1} + 2\cdot 3^{k-2} + 1 < n \leq 3^k$ and $  3^{k} < m. $ Then $A(n,m) = B(n,m).$

\end{prop}

\begin{proof}

Let $(a,b,c)$ be one of the optimal moves from  \cref{maintheorem3c,maintheorem3d,maintheorem3e}.

Then

\begin{align*}
A(n,m) = &  S( a, b, c,  m) 
\\
= &  nm -  A(m,a) - A(m,b) - A(m,c)
\\
= &  nm - \Phi(a) - \Phi(b) - \Phi(c)  
\\
= &  nm - 3^{k-2} \left(1 + 5a - 3^{k-1} - 2 \left\lfloor\frac{a+1}{2}\right\rfloor \right) - 3^{k-2} \left(1 + 5b - 3^{k-1} - 2 \left\lfloor\frac{b+1}{2}\right\rfloor \right) 
\\
&   - 3^{k-2} \left(1 + 5c - 3^{k-1} - 2 \left\lfloor\frac{c+1}{2}\right\rfloor \right) 
\\
= &  nm - 3^{k-2} \left(3 + 5n - 3^{k} - 2 \left\lfloor\frac{a+1}{2}\right\rfloor - 2 \left\lfloor\frac{b+1}{2}\right\rfloor - 2 \left\lfloor\frac{c+1}{2}\right\rfloor \right)
\\
= &  nm - 3^{k-2} \left(3 + 5n - 3^{k} - 2 \left(\left\lfloor\frac{n-1}{2}\right\rfloor + 2\right)  \right)
\\
= &   nm + 2\cdot 3^{k-2}\left\lfloor\frac{n-1}{2}\right\rfloor - 3^{k-2}( 5n -1 -  3^{k} )
\\
= &  B(n,m).\qedhere
\end{align*}
\end{proof}

Combining the above proves $A(n,m) = B(n,m)$ and that $\mathcal{S}$ is an optimal strategy; in other words, we have proved \cref{secthm}.

\end{changemargin}
\end{document}